\documentclass[11pt]{article}
\usepackage{stmaryrd}
\usepackage{mathrsfs}
\usepackage[centertags]{amsmath}
\usepackage{amsfonts, dsfont}
\usepackage{amssymb}
\usepackage{amsthm}

\input xy
\xyoption{all}

\theoremstyle{plain}
\newtheorem{thm}{Theorem}[subsection]
\newtheorem{cor}[thm]{Corollary}
\newtheorem{lem}[thm]{Lemma}
\newtheorem{prop}[thm]{Proposition}

\theoremstyle{definition}
\newtheorem{defn}[thm]{Definition}
\newtheorem{remark}[thm]{Remark}
\newtheorem*{ack}{Acknowledgments}

\numberwithin{equation}{subsection}

\newcommand{\bd}{\begin{defn}}
\newcommand{\ed}{\end{defn}}
\newcommand{\bl}{\begin{lem}}
\newcommand{\el}{\end{lem}}
\newcommand{\bp}{\begin{prop}}
\newcommand{\ep}{\end{prop}}
\newcommand{\bt}{\begin{thm}}
\newcommand{\et}{\end{thm}}
\newcommand{\bc}{\begin{cor}}
\newcommand{\ec}{\end{cor}}
\newcommand{\br}{\begin{remark}}
\newcommand{\er}{\end{remark}}
\newcommand{\bdi}{\begin{diagram}}
\newcommand{\edi}{\end{diagram}}
\newcommand{\beq}{\begin{eqn}}
\newcommand{\eeq}{\end{eqn}}
\newcommand{\ba}{\begin{array}}
\newcommand{\ea}{\end{array}}
\newcommand{\bpf}{\begin{proof}}
\newcommand{\epf}{\end{proof}}

\newcommand{\R}{\mathds{R}}
\newcommand{\Z}{\mathds{Z}}
\newcommand{\Q}{\mathds{Q}}
\newcommand{\Zp}{\mathds{Z}_{p}}
\newcommand{\Qp}{\mathds{Q}_{p}}
\newcommand{\al}{\alpha}
\newcommand{\Ga}{\Gamma}

\newcommand{\La}{\Lambda}
\newcommand{\la}{\lambda}

\newcommand{\Op}{\mathcal{O}}
\newcommand{\ord}{\mathrm{ord}}
\newcommand{\m}{\mathfrak{m}}
\newcommand{\M}{\mathfrak{M}}
\newcommand{\A}{\mathcal{A}}

\DeclareMathOperator{\Sel}{Sel} \DeclareMathOperator{\Gal}{Gal}
\DeclareMathOperator{\Hom}{Hom} \DeclareMathOperator{\rank}{rank}
\DeclareMathOperator{\Ext}{Ext}

\newcommand{\ot}{\otimes}
\newcommand{\ilim}{\displaystyle \mathop{\varinjlim}\limits}
\newcommand{\plim}{\displaystyle \mathop{\varprojlim}\limits}

\newcommand{\coker}{\mathrm{coker}\,}

\newcommand{\cyc}{\mathrm{cyc}}
\newcommand{\cts}{\mathrm{cts}}

\newcommand{\lra}{\longrightarrow}

\newcommand{\ps}[1]{\llbracket #1 \rrbracket}

\begin{document}

\title{Comparing the $\pi$-primary submodules of the dual Selmer
groups}
\author{Meng Fai Lim\footnote{School of Mathematics and Statistics $\&$ Hubei Key Laboratory of Mathematical Sciences,
Central China Normal University, Wuhan, 430079, P.R.China.
 E-mail: \texttt{limmf@mail.ccnu.edu.cn}}}
\date{}
\maketitle

\begin{abstract} \footnotesize

\noindent In this paper, we compare the structure of Selmer groups
of certain classes of Galois representations over an admissible
$p$-adic Lie extension. Namely, we show that the $\pi$-primary
submodules of the Pontryagin dual of the Selmer groups of two Galois
representations have the same elementary representations when two
Galois representations in question are either Tate dual to each
other or are congruent to each other. In the first situation, our
result gives a partial answer to the question of Greenberg on
whether the Pontryagin dual of the Selmer groups of two Galois
representations that are Tate dual to each other are
pseudo-isomorphic (up to a twist of the Iwasawa algebra). In the
second situation, our result will be applied to study the variation
of the $\pi$-primary submodules of the dual Selmer groups of certain
specialization of a big Galois representation. One of the important
ingredient in our proofs is an asymptotic formula for $\pi$-primary
modules over a noncommutative Iwasawa algebra which can be viewed as
a generalization of a weak analog of the classical Iwasawa
asymptotic formula.



\end{abstract}

\section{Introduction}

Throughout the paper, $p$ will always denote a rational prime. Let
$\Op$ be the ring of integers of a fixed finite extension $K$ of
$\Qp$, and let $\pi$ be a local parameter of $\Op$. Let $F$ be a
number field and $F_{\infty}$ an admissible $p$-adic Lie extension
of $F$ with Galois group $G$. For a Galois representation defined
over a number field $F$ with coefficients in $\Op$, one can attach a
Selmer group to these data, and the resulting Selmer group has a natural module
structure over $\Op\ps{G}$ which is the Iwasawa algebra of the
Galois group of the given admissible $p$-adic Lie extension. For a module
over such an Iwasawa algebra, Howson \cite{Ho} and Venjakob
\cite{V02} independently developed the notion of a generalized
$\mu$-invariant which extends the classical notion of the
$\mu$-invariant. Building on this notion of the $\mu$-invariant,
they have established a structural description of the $\pi$-primary
submodules of modules defined over such an Iwasawa algebra (see
\cite{Ho2, V02}). To be precise, they are able to show that the
$\pi$-primary submodules can be expressed uniquely in term of a
product of factors of the form $\Op\ps{G}/\pi^{\al}$ (up to
pseudo-isomorphism). Following the classical situation, we call this
the elementary representation of the $\pi$-primary submodule. In the
midst of establishing their structural theorem, they also show that
the $\mu$-invariant can be expressed as a sum of the $\al$'s
appearing in the elementary representation. The aim of this article
is to compare the $\pi$-primary submodules of the Pontryagin dual of
the Selmer groups of two Galois representations in the two
situations as mentioned in the abstract. We briefly describe these
situations in the next two paragraphs.

The main conjecture of Iwasawa theory is a conjecture on the
relation between a Selmer group and a conjectural $p$-adic
$L$-function (see \cite{CFKSV, FK, G89}). This $p$-adic $L$-function
is expected to satisfy a conjectural functional equation in a
certain sense. In view of the main conjecture and this conjectural
functional equation, one would expect to have certain algebraic
relationship between the Selmer group attached to a Galois
representation and the Selmer group attached to the Tate twist of
the dual of the Galois representation which can be thought as an
algebraic manifestation of the functional equation. In particular,
an immediate consequence of this algebraic functional equation is
that the Selmer group attached to a Galois representation and the
Selmer group attached to the Tate twist of the dual representation
have the same generalized Iwasawa $\mu$-invariants. In the case of a
cyclotomic $\Zp$-extension, this study on the $\mu$-invariants has
been undertaken in \cite{G89, Mat}. (Actually, in \cite{G89},
Greenberg also established the full ``algebraic" functional equation
of the Selmer groups, which we will not treat in this article.
Readers interested in this subject may refer to \cite{BZ, Hs, JP,
LLTT, Z08, Z10}.) For a noncommutative $p$-adic Lie extension, this
equality has been verified for Selmer groups of an abelian variety
without complex multiplication (see \cite{Bh}). In a previous
unpublished note \cite{LimMu}, the author has also verified the equality of
the $\mu$-invariants. In this paper, we will show an even stronger
relation, namely, the $\pi$-primary submodules of the dual Selmer
groups in question have the same elementary representations. In
particular, our result will imply the equality of the generalized
Iwasawa $\mu$-invariants. We like to emphasize that although this
latter result is motivated by the main conjecture of Iwasawa and the
functional equation of the conjectural $p$-adic $L$-functions, we do
not assume these conjectures in all our argument. We also mention
that our result answers the $\pi$-primary part of the question of
Greenberg on whether the Pontryagin dual of the Selmer groups of two
Galois representations are pseudo-isomorphic (up to a twist of the
Iwasawa algebra).

The second situation is concerned with comparing the Selmer groups
of two congruent Galois representations, or in general, a family of
Galois representations with suitable congruence relations. Such
studies were carried out over the cyclotomic $\Zp$-extension in
\cite{EPW, Gr94, GV, Ha, We} and over noncommutative $p$-adic Lie
extensions in \cite{Ch, SS}. One of the motivations behind these
studies lies in the philosophy that the ``Iwasawa main conjecture"
should be preserved by congruences. We should mention that in the
cyclotomic case, this philosophy is rather well understood (see
\cite{EPW, GV, Oc}), although the general noncommutative situation
still seems much of a mystery (for instance, see \cite{B, CS12}). An
important observation made in most of the above cited works is that
if the Iwasawa $\mu$-invariant of one of the Selmer groups vanishes,
so does the other. It is then natural to consider the situation when
the said Iwasawa $\mu$-invariants are nonzero and ask if one can
relate the Iwasawa $\mu$-invariants. To the best of the author's
knowledge, such studies have only been considered in \cite{B, BS}.
As observed in these works and as we will see in this paper, to be
able to even compare the Iwasawa $\mu$-invariants meaningfully, we
require the congruence of the Galois representations to be high
enough. In this paper, we will show the stronger assertion, namely,
that the $\pi$-primary submodules of the Selmer groups of congruent
Galois representations have the same elementary representations.

We now describe briefly the idea behind the comparison of the
$\pi$-primary submodules of the dual Selmer groups of two Galois
representations defined over some number field $F$. Let $A$ and $B$
be quotient modules of the two Galois representations. Here the
quotient module is obtained from the Galois representation by taking
the quotient of the representation by a Galois-invariant lattice.
Then one can attach dual Selmer groups to these quotient modules
(see Section \ref{Arithmetic Preliminaries}) which we denote to be
$X(A/F_{\infty})$ and $X(B/F_{\infty})$, where $F_{\infty}$ is an
admissible $p$-adic Lie extension whose Galois group is a uniform
pro-$p$ group of dimension $r$. It follows from the theory developed
in Subsection \ref{pi-primary modules} that to compare $\pi$-primary
submodules of the dual Selmer groups, we need to show the equality
\[  \mu_{\Op\ps{G}}\Big(X(A/F_{\infty})/\pi^{n}\Big) =
\mu_{\Op\ps{G}}\Big(X(B/F_{\infty})/\pi^{n}\Big) \] for enough $n$
(see Propositions \ref{pseudo-isomorphic2} and
\ref{pseudo-isomorphic3} for details). To prove the above equality,
we proceed by establishing the following estimate

\[ \bigg| \Big(\mu_{\Op\ps{G}}\big(X(A/F_{\infty})/\pi^n\big)-
\mu_{\Op\ps{G}}\big(X(B/F_{\infty})/\pi^n\big)\Big)p^{rm} \bigg| =
O(p^{(r-1)m})  \] for each fixed $n$.  If we denote $F_m$ to be the
intermediate extension of $F$ in $F_{\infty}$ corresponding to fixed
field of the $(m+1)$-term of the lower $p$-series of
$\Gal(F_{\infty}/F)$, we then have the following inequality

\[ \ba{l}
  \bigg| \Big(\mu_{\Op\ps{G}}\big(X(A/F_{\infty})/\pi^n\big)-
\mu_{\Op\ps{G}}\big(X(B/F_{\infty})/\pi^n\big)\Big)p^{rm} \bigg|
\leq \\
 \hspace{1in} \Big| \mu_{\Op\ps{G}}\big(X(A/F_{\infty})/\pi^n\big)p^{rm} -
 \ord_q\big(S(A/F_{\infty})[\pi^n]^{G_m}\big)\Big| \\
 \hspace{1in} +
  \Big| \mu_{\Op\ps{G}}\big(X(B/F_{\infty})/\pi^n\big)p^{rm} -
 \ord_q\big(S(B/F_{\infty})[\pi^n]^{G_m}\big)\Big| \\
 \hspace{1in} + \Big| \ord_q \big(S(A[\pi^n]/F_m)\big) -
 \ord_q\big(S(A/F_{\infty})[\pi^n]^{G_m}\big)\Big| \\
 \hspace{1in} +
 \Big| \ord_q \big(S(B[\pi^n]/F_m)\big) -
 \ord_q\big(S(B/F_{\infty})[\pi^n]^{G_m}\big)\Big| \\
 \hspace{1in} + \Big| \ord_q \big(S(A[\pi^n]/F_m)\big) -
 \ord_q\big(S(B[\pi^n]/F_m)\big)\Big|.
 \ea   \]

 Therefore, we are reduced to showing that the five quantities on
the right are $O(p^{(r-1)m})$ (for a fixed $n$). The estimates for
the first and second quantities are where we require the asymptotic
formulas for $\pi$-primary modules. The estimates for the third and
fourth quantities follow from a descent argument. The estimate of
the final quantity is where we make use of the facts that $A$ and
$B$ are either Tate dual to each other or congruent to each other
for a high enough power. In the case that $A$ and $B$ are Tate dual
to each other, our argument follows the approach of Greenberg in
\cite{G89}. We like to highlight that in the situation of a
cyclotomic $\Zp$-extension, the error quantities have bounded order
in the intermediate sub-extensions of the cyclotomic
$\Zp$-extension. However, since we may have infinite decomposition
of primes over a $p$-adic Lie extension of dimension $r>1$, the
error quantities may not be bounded, and therefore, we will need a
slightly more careful argument.

We now give a brief description of the layout of the paper. In
Section \ref{Algebraic Preliminaries}, we recall certain algebraic
notion which will be used subsequently in the paper. It is also here
where we develop a method to compare the $\pi$-primary submodules of
two $\Op\ps{G}$-modules. In the final subsection, we will prove an
asymptotic formula for $\pi$-primary modules over an Iwasawa algebra
of a (possibly noncommutative) uniform pro-$p$ group. As seen in the
previous paragraph, this asymptotic formula is crucial in
establishing our main results.

In Section \ref{Arithmetic Preliminaries}, we introduce the Selmer
groups which are the main object of study in this paper. These
Selmer groups are defined for a set of data which arises from an
ordinary Galois representation. Actually, to be precise, the Selmer
group that we consider is called the strict Selmer group in
Greenberg's terminology \cite{G89}. We will also introduce another
variant of the Selmer group (called the Greenberg Selmer group) and
an appropriate Selmer complex closely related to the strict Selmer
group we are working with. In Section \ref{pi-submodules of dual
Selmer groups}, we will prove our main theorems for the strict
Selmer groups. But it will follow from the discussion of Section
\ref{Arithmetic Preliminaries} that all our main results also hold
for the Greenberg Selmer group and the Selmer complex. In Section
\ref{Miscellaneous}, we will apply our main results to study the
variation of the $\pi$-primary submodules of the dual Selmer groups
of the specializations of a big Galois representation. When the
$p$-adic Lie extension is of dimension 2, we will also apply our
main results to study the variation of the $\M_H(G)$-property of the
dual Selmer groups.

\begin{ack}
A significant portion of this work was written up when the author
was a Postdoctoral fellow at the GANITA Lab at the University of
Toronto. He would like to acknowledge the hospitality and conducive
working conditions provided by the GANITA Lab and the University of
Toronto while this work was in progress. The author would also like
to thank Otmar Venjakob for answering his question on the asymptotic
formulas. The author is also supported by the National Natural
Science Foundation of China under the Research Fund for
International Young Scientists (Grant No: 11550110172).
        \end{ack}

\section{Algebraic Preliminaries} \label{Algebraic Preliminaries}

In this section, we recall some algebraic preliminaries that will be
required in the later part of the paper. Namely, we gather various
notation and definitions which will be required for the discussion
of the paper. In Subsection \ref{pi-primary modules}, we develop
various criterions which allow us to compare the structure of the
$\pi$-primary submodules of two modules. In the final subsection of
this section, we will prove a asymptotic formula over an Iwasawa
algebra of a uniform pro-$p$ group which is a generalization of a
weak analog to the classical Iwasawa asymptotic formula over an
Iwasawa algebra of $\Ga\cong\Zp$ \cite{Iw}.

\subsection{Compact $p$-adic Lie group} \label{Compact p-adic Lie
group}

Fix a prime $p$. In this subsection, we recall some facts about
compact $p$-adic Lie groups. The standard references for the
material presented here are \cite{DSMS, Laz}.

For a finitely generated pro-$p$ group $G$, we write $G^{p^i} =
\langle g^{p^i}|~g\in G\rangle$, that is, the group generated by the
$p^i$th-powers of elements in $G$. The pro-$p$ group $G$ is said to
be \textit{powerful} if $G/\overline{G^{p}}$ is abelian for odd $p$,
or if $G/\overline{G^{4}}$ is abelian for $p=2$. We define the lower
$p$-series by $P_{1}(G) = G$, and
\[ P_{i+1}(G) = \overline{P_{i}(G)^{p}[P_{i}(G),G]}, ~\mbox{for}~ i\geq 1. \] It follows from \cite[Thm.\
3.6]{DSMS} that if $G$ is a powerful pro-$p$ group, then $G^{p^i} =
P_{i+1}(G)$ and the $p$-power map
\[ P_{i}(G)/P_{i+1}(G)\stackrel{\cdot p}{\lra}
P_{i+1}(G)/P_{i+2}(G)\] is surjective for each $i\geq 1$. If the
$p$-power maps are isomorphisms for all $i\geq 1$, we say that $G$
is \textit{uniformly powerful} (abrev.\ \textit{uniform}). Note that
in this case, we have an equality $|G:P_2(G)| = |P_i(G):
P_{i+1}(G)|$ for every $i\geq 1$. In fact, it is not difficult to
see that $|G: P_{i+1}(G)| = p^{ir}$, where $r= \dim G$.

We now recall the following characterization of compact $p$-adic Lie
groups due to Lazard \cite{Laz} (see also \cite[Cor.\ 8.34]{DSMS}):
a topological group $G$ is a compact $p$-adic Lie group if and only
if $G$ contains a open normal uniform pro-$p$ subgroup. Furthermore,
if $G$ is a compact $p$-adic Lie group without $p$-torsion, it
follows from \cite[Corollaire 1]{Ser} (see also \cite[Chap.\ V
Sect.\ 2.2)]{Laz}) that $G$ has finite $p$-cohomological dimension.

\subsection{Torsion modules and pseudo-null modules} \label{Rank subsection}
As before, $p$ will denote a fixed prime. Let $\Op$ be the ring of
integers of a finite extension of $\Qp$. For a compact $p$-adic Lie
group $G$, the completed group algebra of $G$ over $\Op$ is given by
 \[ \Op\ps{G} = \plim_U \Op[G/U], \]
where $U$ runs over the open normal subgroups of $G$ and the inverse
limit is taken with respect to the canonical projection maps.

When $G$ is pro-$p$ and has no $p$-torsion, it is well known that
$\Op\llbracket G\rrbracket$ is an Auslander regular ring (cf.
\cite[Theorem 3.26]{V02}; for the definition of Auslander regular
rings, see \cite[Definition 3.3]{V02}). Furthermore, the ring
$\Op\ps{G}$ has no zero divisors (cf.\ \cite{Neu}), and therefore,
admits a skew field $Q(G)$ which is flat over $\Op\ps{G}$ (see
\cite[Chapters 6 and 10]{GW} or \cite[Chapter 4, \S 9 and \S
10]{Lam}). If $M$ is a finitely generated $\Op\ps{G}$-module, we
define the $\Op\ps{G}$-rank of $M$ to be
$$ \rank_{\Op\ps{G}}M  = \dim_{Q(G)} Q(G)\ot_{\Op\ps{G}}M. $$
 We will say that an $\Op\ps{G}$-module $M$ is \textit{torsion} if
$\rank_{\Op\ps{G}} M = 0$. As we will also need to work with various
equivalent formulations of a torsion $\Op\ps{G}$-module, we state
the following.

\bl \label{torsion is Hom zero} Let $\La$ be a Auslander regular
ring with no zero divisors. Let $M$ be a finitely generated
$\La$-module. Then the following are equivalent.

\smallskip $(a)$ The canonical map
 $\phi: M\lra \Hom_{\La}(\Hom_{\La}(M,\La), \La)$ is zero.

\smallskip $(b)$ $Q(\La)\ot_{\La}M =0$, where $Q(\La)$ is the skew
field of $\La$.

\smallskip $(c)$ $\Hom_{\La}(M,\La)= 0$.
 \el

\bpf The equivalence of (a) and (c) follows from \cite[Remark
3.7]{V02}. Suppose that $Q(\La)\ot_{\La}M = 0$. Let $f\in
\Hom_{\La}(M,\La)$ and $x\in M$. Then since $Q(\La)\ot_{\La}M = 0$,
there exists $\la\in\La\setminus\{0\}$ such that $\la x=0$. This in
turn implies that $\la f(x) = f(\la x) = 0$. Since $\La$ has no zero
divisor, we have $f(x)=0$. This shows that $\Hom_{\La}(M,\La)= 0$
and the implication (b)$\Rightarrow$(c). Conversely, suppose that
$\Hom_{\La}(M,\La)= 0$. Write $M^{++} =
\Hom_{\La}(\Hom_{\La}(M,\La), \La)$. By \cite[Proposition 2.5]{V02}
and the Auslander condition, the canonical map
 $\phi: M\lra M^{++}$
has kernel and cokernel which are $\La$-torsion. Therefore, $\phi$
induces an isomorphism
$$Q(\La)\ot_{\La}M \stackrel{\sim}{\lra} Q(\La)\ot_{\La}M^{++}.$$
Now if $\phi =0$, then it will follow that $Q(\La)\ot_{\La}M = 0$.
This establishes (a)$\Rightarrow$(b). \epf

Therefore, if $G$ is pro-$p$ and has no $p$-torsion, it follows from
the above lemma that a finitely generated $\Op\ps{G}$-module $M$ is
torsion if and only if $\Hom_{\Op\ps{G}}(M, \Op\ps{G}) =0$. Now if
$M$ is a finitely generated torsion $\Op\ps{G}$-module, we say that
$M$ is \textit{pseudo-null} if $\Ext^1_{\Op\ps{G}}(M, \Op\ps{G})
=0$. For an equivalent definition, we refer readers to
\cite[Definitions 3.1 and 3.3; Proposition 3.5(ii)]{V02}. For the
purpose of this article, the definition we adopt will suffice.
Finally, we mention that every subquotient of a torsion
$\Op\ps{G}$-module (resp., pseudo-null $\Op\ps{G}$-module) is also
torsion (resp. pseudo-null).

\subsection{$\mu$-invariant}

Let $\Op$ be the ring of integers of a fixed finite extension $K$ of
$\Qp$ as defined in the preceding subsection. Fix a local parameter
$\pi$ for $\Op$ and denote the residue field of $\Op$ by $k$. The
completed group algebra of $G$ over $k$ is given by
 \[ k\ps{G} = \plim_U k[G/U], \]
where $U$ runs over the open normal subgroups of $G$ and the inverse
limit is taken with respect to the canonical projection maps.

For a compact $p$-adic Lie group $G$ without $p$-torsion, it follows
from \cite[Theorem 3.30(ii)]{V02} that $k\ps{G}$ is an Auslander
regular ring. Furthermore, if $G$ is pro-$p$ without $p$-torsion,
then the ring $k\ps{G}$ has no zero divisors (cf. \cite[Theorem
C]{AB}). Therefore, one can define the notion of $k\ps{G}$-rank as
above when $G$ is pro-$p$ without $p$-torsion. We will say that that
the module $N$ is a \textit{torsion} $k\ps{G}$-module if
$\rank_{k\ps{G}}N = 0$. By Lemma \ref{torsion is Hom zero}, we have
that $N$ is a torsion $k\ps{G}$-module if and only if
$\Hom_{k\ps{G}}(N, k\ps{G}) = 0$.

For a given finitely generated $\Op\ps{G}$-module $M$, we denote
$M(\pi)$ to be the $\Op\ps{G}$-submodule of $M$ which consists of
elements of $M$ that are annihilated by some power of $\pi$. Since
the ring $\Op\ps{G}$ is Noetherian, the module $M(\pi)$ is finitely
generated over $\Op\ps{G}$. Therefore, one can find an integer
$r\geq 0$ such that $\pi^r$ annihilates $M(\pi)$. Following
\cite[Formula (33)]{Ho}, we define
  \[\mu_{\Op\ps{G}}(M) = \sum_{i\geq 0}\rank_{k\ps{G}}\big(\pi^i
   M(\pi)/\pi^{i+1}\big). \]
(For another alternative, but equivalent, definition, see
\cite[Definition 3.32]{V02}.) By the above discussion and our
definition of $k\ps{G}$-rank, the sum on the right is a finite one.
It is clear from the definition that $\mu_{\Op\ps{G}}(M) =
\mu_{\Op\ps{G}}(M(\pi))$. Also, it is not difficult to see that this
definition coincides with the classical notion of the
$\mu$-invariant for $\Ga$-modules when $G=\Ga\cong\Zp$.

We now record certain properties of the $\mu_{\Op\ps{G}}$-invariant
which will be required in the subsequent of the paper.

\bl \label{mu lemma} Let $G$ be a compact pro-$p$ $p$-adic Lie group
with no $p$-torsion and let $M$ be a finitely generated
$\Op\ps{G}$-module. Then we have the following statements.

\begin{enumerate}
\item[$(a)$] For every finitely generated $\Op\ps{G}$-module $M$, one has
$$\mu_G(M) = \sum_{i\geq 0} (-1)^i\mathrm{ord}_q\big(H_i(G,M(\pi))\big)£¬$$
 where $q$ is the cardinality of $k$.

\item[$(b)$]  Suppose that $G$ has a closed normal subgroup $H$ such that
$G/H\cong \Zp$. If $M$ is a $\Op\ps{G}$-module which is finitely
generated over $\Op\ps{H}$, then one has $\mu_{\Op\ps{G}}(M) =0$.

\item[$(c)$] Suppose that we are given a short exact sequence of finitely generated
$\Op\ps{G}$ modules
\[ 0\lra M'\lra M\lra M'' \lra 0.\]

\begin{enumerate}
\item[$(1)$] One has $\mu_{\Op\ps{G}}(M) \leq \mu_{\Op\ps{G}}(M') + \mu_{\Op\ps{G}}(M'')$.
Moreover, if $M$, and hence also $M'$ and $M''$, is
$\Op\ps{G}$-torsion, the inequality is an equality.

\item[$(2)$] If $\mu_{\Op\ps{G}}(M'')=0$, then one has $\mu_{\Op\ps{G}}(M') =
\mu_{\Op\ps{G}}(M)$.
\end{enumerate}

\item[$(d)$]  Suppose that $G$ has a closed normal subgroup $H$ such that
$G/H\cong \Zp$ and suppose that we are given an exact sequence of
finitely generated $\Op\ps{G}$-modules
\[ A\lra B \lra C\lra D\]
such that $A$ is finitely generated over $\Op\ps{H}$ and
$\mu_G(D)=0$. Then one has the equality $\mu_G(B) = \mu_G(C)$.

\item[$(e)$]  $\mu_{\Op\ps{G}}(M) =
0$ if and only if $M(\pi)$ is pseudo-null.
  \end{enumerate}
\el

\bpf Statements (a), (b) and (c)(1) are proven in \cite[Corollary
1.7]{Ho}, \cite[Lemma 2.7]{Ho} and \cite[Proposition 1.8]{Ho}
respectively. Statement (e) is shown in \cite[Remark 3.33]{V02}. The
remaining statements can be deduced from the previous statements
without too much difficulties. \epf

\subsection{$\pi$-primary modules} \label{pi-primary modules}

Throughout this subsection, $G$ will always denote a pro-$p$ $p$-adic Lie
group without $p$-torsion. Therefore, both rings $\Op\ps{G}$ and $k\ps{G}$ are Auslander
regular and have no zero divisors. For a finitely generated
$\Op\ps{G}$-module $M$, it then follows from \cite[Proposition
1.11]{Ho2} (see also \cite[Theorem 3.40]{V02}) that there is a
$\Op\ps{G}$-homomorphism
\[ \varphi: M(\pi) \lra \bigoplus_{i=1}^s\Op\ps{G}/\pi^{\al_i},\] whose
kernel and cokernel are pseudo-null $\Op\ps{G}$-modules, and where
the integers $s$ and $\al_i$ are uniquely determined. We will call
$\bigoplus_{i=1}^s\Op\ps{G}/\pi^{\al_i}$ the \textit{elementary
representation} of $M(\pi)$. In fact, in the process of establishing
the above, one also has that $\mu_{\Op\ps{G}}(M) = \displaystyle
\sum_{i=1}^s\al_i$ (see loc. cit.). We will set
\[\theta_{\Op\ps{G}}(M) : = \max_{1\leq i\leq s}\{\al_i\}.\]

The following fundamental lemma gives a relationship between the
$\mu_{\Op\ps{G}}$-invariant and $\Op\ps{G}$-rank of a finitely
generated $\Op\ps{G}$-module.

\bl \label{mu and rank} Let $M$ be a finitely generated
$\Op\ps{G}$-module. Suppose that there is a $\Op\ps{G}$-homomorphism
\[ \varphi: M(\pi) \lra \bigoplus_{i=1}^s\Op\ps{G}/\pi^{\al_i},\] whose
kernel and cokernel are pseudo-null $\Op\ps{G}$-modules. Then we
have
\[\mu_{\Op\ps{G}}(M/\pi^n) =
n\rank_{\Op\ps{G}}(M) + \sum_{i=1}^s\min\{n,\al_i\} \quad
\mbox{for}~ n\geq 1.
\] In particular, we have $\mu_{\Op\ps{G}}(M/\pi^n) \leq
n\rank_{\Op\ps{G}}(M)+ \mu_{\Op\ps{G}}(M)$ which is an equality if
and only if $n\geq \theta_{\Op\ps{G}}(M)$. \el

\bpf
  Write $M_f = M/ M(\pi)$. Consider the following commutative
  diagram
  \[ \SelectTips{eu}{}
\xymatrix{
  0 \ar[r] & M(\pi) \ar[d]^{\pi^n} \ar[r]^{} & M \ar[d]^{\pi^n}
  \ar[r]^{} & M_f \ar[d]^{\pi^n} \ar[r]& 0 \\
  0 \ar[r] & M(\pi) \ar[r] & M
  \ar[r] & M_f \ar[r] &0   }
\] with exact rows, and the vertical maps are given by
multiplication by $\pi^n$. Since $M_f$ has no $\pi$-torsion, the
rightmost vertical map is injective, and therefore, we have an exact
sequence
\[ 0\lra M(\pi)/\pi^n\lra M/\pi^n \lra M_f/\pi^n \lra 0 \]
of torsion $\Op\ps{G}$-modules. By Lemma \ref{mu lemma}(c)(1), we
have
 \[ \mu_{\Op\ps{G}}(M/\pi^n) = \mu_{\Op\ps{G}}(M(\pi)/\pi^n) +
 \mu_{\Op\ps{G}}(M_f/\pi^n).\]
To prove the lemma, it therefore suffices to show the following two
equalities.

(1) $\mu_{\Op\ps{G}}(M(\pi)/\pi^n) = \sum_{i=1}^s\min\{n,\al_i\}$.

(2) $\mu_{\Op\ps{G}}(M_f/\pi^n) = n\rank_{\Op\ps{G}}(M_f)
=n\rank_{\Op\ps{G}}(M)$.

To see that (1) holds, note that since any subquotient of a
pseudo-null module is also pseudo-null, it follows that $\varphi$
induces an $\Op\ps{G}$-homomorphism
\[ M(\pi)/\pi^n \lra \bigoplus_{i=1}^s\Op\ps{G}/\pi^{\min\{n,\al_i\}},\] whose
kernel and cokernel are pseudo-null $\Op\ps{G}$-modules. The
equality in (1) will now follow by combining this observation with
statements (c)(1) and (e) of Lemma \ref{mu lemma}.

Since $M(\pi)$ is clearly a torsion $\Op\ps{G}$-module, we have
$\rank_{\Op\ps{G}}(M_f) = \rank_{\Op\ps{G}}(M)$. Therefore, it
remains to verify the first equality in (2). In other words, we are
reduced to showing that if $M$ is a finitely generated
$\Op\ps{G}$-module with $M(\pi) = 0$, then $\mu_{\Op\ps{G}}(M/\pi^n)
= n\rank_{\Op\ps{G}}(M)$. We shall proceed by induction. Suppose
that $n=1$. Then we have
\[\rank_{\Op\ps{G}}(M) =
\rank_{k\ps{G}}(M/\pi) = \sum_{i\geq 0}
(-1)^i\dim_k\big(H_i(G,M/\pi)\big) = \mu_{\Op\ps{G}}(M/\pi), \]
where the first equality follows from \cite[Corollary 1.10]{Ho} and
the assumption that $M[\pi] = 0$, the second equality follows from
\cite[Proposition 1.6]{Ho}, and the third equality follows from
Lemma \ref{mu lemma}(a). Therefore, we have established the $n=1$
case.

Now suppose that $n>1$, and suppose that
$\mu_{\Op\ps{G}}(M/\pi^{n-1}) = (n-1)\rank_{\Op\ps{G}}(M)$. Then
consider the following commutative diagram
 \[ \SelectTips{eu}{}
\xymatrix{
   & M \ar[d]^{\pi^{n-1}} \ar@{=}[r] & M \ar[d]^{\pi^n}
   &  &  \\
  0 \ar[r] & M \ar[r]^{\pi} & M
  \ar[r] & M/\pi \ar[r] &0   } \]
with exact bottom row (note that injectivity follows from the
assumption that $M[\pi]= 0$). By the snake lemma, we have an exact
seqeuence
\[ 0\lra M/\pi^{n-1} \lra M/\pi^n \lra M/\pi \lra 0 \]
of torsion $\Op\ps{G}$-modules which in turn yields
 \[\ba{rl}
  \mu_{\Op\ps{G}}(M/\pi^n) ~= &\! \mu_{\Op\ps{G}}(M/\pi) +
 \mu_{\Op\ps{G}}(M/\pi^{n-1}) \\
 ~= &\! \rank_{\Op\ps{G}}(M) + (n-1)\rank_{\Op\ps{G}}(M) ~=~
 n\rank_{\Op\ps{G}}(M).
  \ea \] The proof of the lemma is completed. \epf

\br
 When $G\cong\Zp^r$, one can prove the above lemma by appealing
 directly to the structure theory (cf. \cite[Proposition 5.1.7]{NSW}).
 \er

We record the following lemma which enables one to relate $M(\pi)$
and $N(\pi)$ in certain situation.

\bl \label{pseudo-isomorphism lemma}
 Suppose that $H$ is a closed normal subgroup of $G$ with $G/H\cong \Zp$.
Let $\varphi: M \lra N$ be a homomorphism of finitely generated
$\Op\ps{G}$-modules, whose kernel and cokernel are finitely
generated over $\Op\ps{H}$. Then $M(\pi)$ and $N(\pi)$ have the same
elementary representations. \el

\bpf
  The statement will follow if it holds in the two special cases of exact
sequences
\[\ba{c} 0\lra P \lra M \lra N\lra 0, \\
  0\lra M \lra N \lra P\lra 0, \ea
\] where $P$ is a finitely generated $\Op\ps{H}$-module. We will
prove the first case, the second case has a similar argument. Choose
a sufficiently large $n$ such that $\pi^n$ annihilates $M(\pi)$ and
$N(\pi)$. Consider the following commutative diagram
  \[ \SelectTips{eu}{}
\xymatrix{
  0 \ar[r] & P \ar[d]^{\pi^n} \ar[r]^{} & M \ar[d]^{\pi^n}
  \ar[r]^{} & N \ar[d]^{\pi^n} \ar[r]& 0 \\
  0 \ar[r] & P \ar[r] & M
  \ar[r] & N \ar[r] &0   }
\] with exact rows, and the vertical maps are given by
multiplication by $\pi^n$. Applying the Snake Lemma, we obtain
\[ 0\lra P[\pi^n]\lra M(\pi) \lra N(\pi) \lra P/\pi^n \]
By Lemma \ref{mu lemma}(b) and (e), we have that $P[\pi^n]$ and
$P/\pi^n$ are pseudo-null $\Op\ps{G}$-modules. Let
\[ f: N(\pi) \lra \bigoplus_{i=1}^s\Op\ps{G}/\pi^{\al_i}\]
 be a homomorphism of $\Op\ps{G}$-modules,
whose kernel and cokernel are pseudo-null $\Op\ps{G}$-modules. Then
 \[f\circ \varphi: M(\pi) \lra
\bigoplus_{i=1}^s\Op\ps{G}/\pi^{\al_i}\] is a homomorphism of
$\Op\ps{G}$-modules, whose kernel and cokernel are pseudo-null.
Therefore, $M(\pi)$ and $N(\pi)$ have the same elementary
representations.
 \epf

\bp \label{mu inequality}
 Let $M$ and $N$ be two finitely generated torsion $\Op\ps{G}$-modules such
that $\mu_{\Op\ps{G}}(M/\pi^{\theta_{\Op\ps{G}}(M)}) =
\mu_{\Op\ps{G}}(N/\pi^{\theta_{\Op\ps{G}}(M)})$. Then we have
\[\mu_{\Op\ps{G}}(M) \leq \mu_{\Op\ps{G}}(N).\] \ep

\bpf
 By Lemma \ref{mu and rank}, we have
 \[ \mu_{\Op\ps{G}}(M) = \mu_{\Op\ps{G}}(M/\pi^{\theta_{\Op\ps{G}}(M)}) =
\mu_{\Op\ps{G}}(N/\pi^{\theta_{\Op\ps{G}}(M)}) \leq
\mu_{\Op\ps{G}}(N). \] \epf

Now for a given finitely generated torsion $\Op\ps{G}$-module $M$,
the elementary representation of $M(\pi)$ can be rewritten as
 \[\bigoplus_{i=1}^{\theta}(\Op\ps{G}/\pi^{i})^{s_i}\]
 for some nonnegative integers $s_i$. Here $\theta =
 \theta_{\Op\ps{G}}(M)$. Then for every
$1\leq n \leq \theta$, we have
\[ \mu_{\Op\ps{G}}(M/\pi^n) = s_1 + 2 s_2 + \cdots + (n-1)s_{n-1} + n(s_n + \cdots
s_{\theta}). \] Putting these equations into a matrix form, we have

\[
\begin{pmatrix}
  \mu_{\Op\ps{G}}(M/\pi) \\
  \mu_{\Op\ps{G}}(M/\pi^2)\\
  \mu_{\Op\ps{G}}(M/\pi^3)\\
  \vdots \\
  \mu_{\Op\ps{G}}(M/\pi^{\theta})
 \end{pmatrix}
=
 \begin{pmatrix}
  1 & 1 & 1& \cdots &1 \\
  1 & 2 & 2& \cdots &2\\
  1 & 2 & 3& \cdots &3\\
  \vdots  & \vdots & \vdots & \ddots & \vdots\\
  1 & 2 & 3& \cdots & \theta
 \end{pmatrix}
  \begin{pmatrix}
  s_1 \\
  s_2\\
  s_3\\
  \vdots \\
  s_{\theta} \end{pmatrix}.
  \]
It is a simple linear algebra exercise to verify that the square
matrix in the above equation is invertible. Therefore, the integers
$s_i$, and hence the elementary representation of $M(\pi)$, are
determined by the values of $\mu_{\Op\ps{G}}(M/\pi^i)$. We record
this observation in the next proposition.

\bp \label{pseudo-isomorphic}
 Let $M$ and $N$ be two finitely generated torsion $\Op\ps{G}$-modules. Then the following
 are equivalent.
 \begin{enumerate}
\item[$(a)$]   $\theta_{\Op\ps{G}}(M)= \theta_{\Op\ps{G}}(N)$, and
$\mu_{\Op\ps{G}}(M/\pi^i) = \mu_{\Op\ps{G}}(N/\pi^i)$ for every
$1\leq i \leq \theta_{\Op\ps{G}}(M)$.
\item[$(b)$]  $M(\pi)$ and $N(\pi)$ have the same elementary representations.
\end{enumerate}\ep

\bpf
 The discussion before the proposition establishes the implication
 $(a)\Rightarrow(b)$. The reverse implication is obvious.
\epf

The preceding proposition may be difficult to apply due to the
condition $\theta_{\Op\ps{G}}(M) = \theta_{\Op\ps{G}}(N)$ which is
perhaps not easy to check. However, one can build on the proposition
to obtain the following which is perhaps easier for application.

\bp \label{pseudo-isomorphic2}
 Let $M$ and $N$ be two finitely generated $\Op\ps{G}$-modules such
that $M$ is a torsion $\Op\ps{G}$-module and such that
$\mu_{\Op\ps{G}}(M/\pi^i) = \mu_{\Op\ps{G}}(N/\pi^i)$ for every
$1\leq i \leq \theta_{\Op\ps{G}}(M)+1$.

Then $N$ is torsion over $\Op\ps{G}$ and we have the equality
$\theta_{\Op\ps{G}}(M) = \theta_{\Op\ps{G}}(N)$. In particular,
$M(\pi)$ and $N(\pi)$ have the same elementary representations. \ep

\bpf We first prove the proposition for the case when
$\theta_{\Op\ps{G}}(M) =0$. Then by the hypothesis of the
proposition, we have $\mu_{\Op\ps{G}}(N/\pi) =
\mu_{\Op\ps{G}}(M/\pi) = 0$. By Lemma \ref{mu and rank}, this in
turns implies that $\rank_{\Op\ps{G}}(N)=0$ and $\mu_{\Op\ps{G}}(N)
=0$. Therefore, we have that $N$ is torsion over $\Op\ps{G}$ and
$\theta_{\Op\ps{G}}(N) =0$. Hence we have that $M(\pi)$ and $N(\pi)$
are both pseudo-null by Lemma \ref{mu lemma}(e), and therefore, have
the same elementary representations.

Now suppose that $\theta_{\Op\ps{G}}(M) \geq 1$. Suppose that the elementary factor
of $N(\pi)$ is given by
\[\bigoplus_{i=1}^t\Op\ps{G}/\pi^{\beta_i}.\]
 Write $a= \rank_{\Op\ps{G}}(N)$. By Lemma \ref{mu and rank}, we then have
 \[ \mu_{\Op\ps{G}}(M) = \mu_{\Op\ps{G}}(M/\pi^n) =
\mu_{\Op\ps{G}}(N/\pi^{n}) = na+ \sum_{i=1}^t\min\{n ,\beta_i\}
\] for $n = \theta_{\Op\ps{G}}(M), \theta_{\Op\ps{G}}(M)+1$.
 This in turn implies that
\[\theta_{\Op\ps{G}}(M)a+\sum_{i=1}^t\min\{\theta_{\Op\ps{G}}(M),\beta_i\} =
\big(\theta_{\Op\ps{G}}(M)+1\big)a+\sum_{i=1}^t\min\{\theta_{\Op\ps{G}}(M)+1,\beta_i\}.
\]
 Since one always has $\theta_{\Op\ps{G}}(M)a\leq
\big(\theta_{\Op\ps{G}}(M)+1\big)a$ and
$\min\{\theta_{\Op\ps{G}}(M),\beta_i\} \leq
\min\{\theta_{\Op\ps{G}}(M)+1,\beta_i\}$, in order for the above
equality to hold, we must have $a=0$ and
$\min\{\theta_{\Op\ps{G}}(M),\beta_i\} =
\min\{\theta_{\Op\ps{G}}(M)+1,\beta_i\}$ for $1\leq i\leq t$. The
formal equality then shows that $N$ is a torsion $\Op\ps{G}$-module,
and the latter equalities show that $\beta_i\leq
\theta_{\Op\ps{G}}(M)$ for all $i$, or in other words,
$\theta_{\Op\ps{G}}(N)\leq \theta_{\Op\ps{G}}(M)$. Therefore, we may
repeat the above argument (noting that we have shown that $N$ is
$\Op\ps{G}$-torsion) replacing $\theta_{\Op\ps{G}}(M)$ by
$\theta_{\Op\ps{G}}(N)$ and interchanging the roles of $M$ and $N$
to obtain the reverse inequality $\theta_{\Op\ps{G}}(M)\leq
\theta_{\Op\ps{G}}(N)$. The remaining assertion will now follow from
an application of Proposition \ref{pseudo-isomorphic}.
 \epf

In particular, it follows from Proposition \ref{pseudo-isomorphic2}
that if $M$ and $N$ are two finitely generated torsion
$\Op\ps{G}$-modules such that $\mu_{\Op\ps{G}}(M/\pi^i) =
\mu_{\Op\ps{G}}(N/\pi^i)$ for \textit{every} $i \geq 1$, then
$M(\pi)$ and $N(\pi)$ have the same elementary representations. In
fact, we can even establish the following stronger statement.

\bp \label{pseudo-isomorphic3}
 Let $M$ and $N$ be two finitely generated $\Op\ps{G}$-modules such
that $\mu_{\Op\ps{G}}(M/\pi^i) = \mu_{\Op\ps{G}}(N/\pi^i)$ for every
$i \geq 1$.

Then we have that $\rank_{\Op\ps{G}}(M) = \rank_{\Op\ps{G}}(N)$ and
that $M(\pi)$ and $N(\pi)$ have the same elementary representations.
\ep

\bpf
 For $n\geq \max\{\theta_{\Op\ps{G}}(M), \theta_{\Op\ps{G}}(N)\}$,
it follows from the assumption $\mu_{\Op\ps{G}}(M/\pi^n) =
\mu_{\Op\ps{G}}(N/\pi^n)$ and Lemma \ref{mu and rank} that
\[ n\rank_{\Op\ps{G}}(M) + \mu_{\Op\ps{G}}(M) =
n\rank_{\Op\ps{G}}(N) + \mu_{\Op\ps{G}}(N).\]
 In other words, we have
\[ \rank_{\Op\ps{G}}(M) + \frac{1}{n}\mu_{\Op\ps{G}}(M) =
\rank_{\Op\ps{G}}(N) + \frac{1}{n}\mu_{\Op\ps{G}}(N).\]
 Letting $n\rightarrow\infty$, we obtain
 \[\rank_{\Op\ps{G}}(M) = \rank_{\Op\ps{G}}(N). \] This proves
the first assertion.

As seen in the proof of Lemma \ref{mu and rank}, we have
\[ \mu_{\Op\ps{G}}(M/\pi^i) = \mu_{\Op\ps{G}}(M(\pi)/\pi^i) +
 i \rank_{\Op\ps{G}}(M).\]
 One has a similar equality for $N$. It then follows from these equalities and
what we proved in the preceding paragraph that
\[ \mu_{\Op\ps{G}}(M(\pi)/\pi^i) = \mu_{\Op\ps{G}}(N(\pi)/\pi^i)\]
for all $i\geq 1$. The second assertion will now follow from an
application of Proposition \ref{pseudo-isomorphic2} on $M(\pi)$ and
$N(\pi)$.
 \epf

\subsection{An asymptotic formula}
 In this subsection, $G$ will always denote a pro-$p$ $p$-adic Lie group without $p$-torsion.
We denote by $r$ the dimension of $G$. We fix an open normal uniform subgroup $G_0$ of $G$ (such a group exists by virtue of Lazard's theorem \cite{Laz}). In the event that $G$ is already a uniform group, we take $G_0 = G$. We now write $G_m$ for
$P_{m+1}(G_0)$ which is defined as in Subsection \ref{Compact p-adic
Lie group}. As before, $\Op$ is the ring of integers of a finite
extension $K$ of $\Qp$, $\pi$ is a local parameter of $\Op$ and $k$
is the residue field of $\Op$. Denote $q$ to the order of $k$. Every
finite $\Op$-module can be viewed as a $\Op/\pi^n$-module for some
$n$. Since $\Op/\pi^n$ has order of a power of $q$, so is every
finite $\Op$-module. For a finite $\Op$-module, we will denote
$\ord_q(M)$ to be the exponent of $q$ in the order of $M$, i.e.,
$|M| = q^{\ord_q(M)}$.

We take this opportunity to introduce a notion which will used in
this paper. A sequence of real numbers $(a_m)_{m\geq 1}$ is said to
satisfy $O(Q^m)$ for some nonnegative number $Q$ if $|a_m| \leq
CQ^m$ for some constant $C$ (independent of $m$) for all
sufficiently large $m$. We will write $a_m = O(Q^m)$. If
$(b_m)_{m\geq 1}$ is another sequence of real numbers, we will write
$a_m = b_m + O(Q^m)$ to mean $a_m - b_m = O(Q^m)$.

We can now state the main theorem of this subsection.

\bt \label{asymptotic formula}
 Let $G$ be pro-$p$ $p$-adic Lie group without $p$-torsion. Write $r=\dim G$. Let $M$ be a finitely
generated $\Op\ps{G}$-module such that $M=M(\pi)$. Then we have
 \[ \ord_q \big(M_{G_m}\big) = [G:G_0]\mu_{\Op\ps{G}}(M)p^{rm} + O(p^{(r-1)m}) \]
 and
 \[ \ord_q\big(H_i(G_m, M)\big) = O(p^{(r-1)m}) \] for every $i\geq 1$. \et

\br The first asymptotic formula in the above result is a weak
analog of the asymptotic formula of Iwasawa \cite[Theorem 4]{Iw}
(see also \cite[Proposition 5.3.17]{NSW}). When $G=\Zp^r$, this can
also be viewed as a weak analog of the asymptotic formula of Cucuo
and Monsky \cite[Theorem 4.13]{CM} (see also \cite[Theorem
3.12]{Mon}). \er

It seems possible that the above formula might be known among the
experts. Despite so, due to a lack of proper reference, we will
include a proof here. (In fact, as we shall see, the tools required
for the proof are available from \cite{Har, Har2, Ho,Ho2, V02}.) For
our purpose in this paper, we will only require the first asymptotic
formula. Despite so, we have included the proof of the asymptotic
formulas for the higher cohomology groups for completeness.

The proof of Theorem \ref{asymptotic formula} will take up the
remainder of this subsection. As a start, we note that
$[G:G_0]\mu_{\Op\ps{G}}(M)= \mu_{\Op\ps{G_0}}(M)$. Since a
$\Op\ps{G}$-module can be viewed as a $\Op\ps{G_0}$ by restriction
of scalars, it suffices to prove the theorem under the assumption
that $G$ is uniform. \textit{In view of this, we therefore can and
do assume that $G$ is uniform for the subsequent of this
subsection.} For the preparation of the proof, we require a few
lemmas.

\bl \label{asymptotic G_m}
 Let $M$ be a finitely generated torsion $k\ps{G}$-module. Then
 \[ \ord_q \big(M_{G_m}\big) =  O(p^{(r-1)m}).\]
\el

\bpf
 Since $M$ is finitely generated torsion over $k\ps{G}$,
there is a surjective map
 \[\bigoplus_j k\ps{G}/k\ps{G}f_j\lra M \] for a finite set of
non-zero and non-unital elements $f_i\in k\ps{G}$. Therefore, one is
reduced to the case $M=k\ps{G}/k\ps{G}f$ for some non-zero and
non-unital $f$. The remainder of the proof then proceeds as in the
proofs of \cite[Lemma 1.10.1]{Har} and \cite[Theorem 1.10]{Har2},
where one passes to the graded ring of $k\ps{G}$ and appeals to the
theory of Hilbert polynomials. \epf

We will also require an estimate for $\ord_q \big(H_i(G_m, M)\big)$.
Before showing this, we need the following lemma.

\bl
 Let $f$ be a
 nonzero nonunital element of $k\ps{G}$, and set $M = k\ps{G}/k\ps{G}f$.
 Then for every $m$, we have
 \[ \ord_q\big(H_1(G_m, M)) = \ord_q\big(M_{G_m} \big)\]
  and $H_i(G_m, M)=0$ for $i\geq 2$.
\el

\bpf
 Since $k\ps{G}$ has no zero divisors, we have an exact sequence
 \[ 0\lra k\ps{G} \stackrel{\cdot f}{\lra} k\ps{G} \lra M \lra 0. \]
 Since $H_i(G_m, k\ps{G}) = 0$ for $i\geq
1$, it follows from considering the $G_m$-homology that we obtain an
exact sequence
\[ 0\lra H_1(G_m, M)\lra k[G/G_m] \lra k[G/G_m] \lra M_{G_m} \lra 0, \]
and the vanishing of $H_i(G_m, M)$ for $i\geq 2$. The first
conclusion of the lemma is now immediate from the four term exact
sequence. \epf

We can now give an estimate for $\ord_q \big(H_i(G_m, M)\big)$.

\bl \label{asymptotic H_1 G_m}
 Let $M$ be a finitely generated torsion $k\ps{G}$-module. Then for
 every $i\geq 1$, we have
  \[ \ord_q \big(H_i(G_m, M)\big) =  O(p^{(r-1)m}), \]
  where $r$ denotes the dimension of $G$. \el

\bpf
  As above, we have a exact sequence
 \[ 0\lra N \lra \bigoplus_j k\ps{G}/k\ps{G}f_j\lra M \lra 0 \]
 of torsion $k\ps{G}$-modules. Taking the $G_m$-homology, we have an
 exact sequence
 \[ H_i\Big(G_m, \bigoplus_j k\ps{G}/k\ps{G}f_j\Big) \lra H_i(G_m, M)
 \lra H_{i-1}(G_m, N). \]
The required estimates will follow from the previous two lemmas.
\epf

We now establish our estimates for the case when $M = M(\pi)$ is a
finitely generated pseudo-null $\Op\ps{G}$-module. Note that the
said module has trivial $\mu_{\Op\ps{G}}$-invariant by Lemma \ref{mu
lemma}(e).

\bl \label{asymptotic formula lemma}
 Let $G$ be a uniform pro-$p$ group of dimension $r$.
 Suppose that $M = M(\pi)$ is a finitely generated pseudo-null $\Op\ps{G}$-module.
 Then for each $i\geq 0$, we have
 \[ \ord_q \big(H_i(G_m, M)\big) =  O(p^{(r-1)m}). \]
\el

\bpf Since $M$ is annihilated by a power of $\pi$, it is a finite
successive extension of subquotients $\pi^iM/\pi^{i+1}$. Therefore,
it suffices to bound each of these subquotients. Since subquotients
of a pseudo-null $\Op\ps{G}$-module are also pseudo-null, we are
reduced to showing that if $M$ is a finitely generated pseudo-null
$\Op\ps{G}$-module with $\pi M =0$, then $\ord_q \big(M_{G_m}\big) =
O(p^{(r-1)m})$. Since $\pi M =0$, we may also view $M$ as a
$k\ps{G}$-module. By a standard spectral sequence argument (for
instance, see \cite[Section 3.4]{V02}), we have
\[ \Ext^i_{k\ps{G}}(M, k\ps{G})\cong \Ext^{i+1}_{\Op\ps{G}}(M, \Op\ps{G})\]
for any integer $i$. In particular, we have
\[ \Hom_{k\ps{G}}(M, k\ps{G})\cong \Ext^{1}_{\Op\ps{G}}(M, \Op\ps{G}) = 0,\]
where the last equality follows from the fact that $M$ is
pseudo-null over $\Op\ps{G}$. Hence $M$ is a torsion
$k\ps{G}$-module. The first estimate then follows from Lemma
\ref{asymptotic G_m}. The estimates for the higher cohomology groups
can be proven similarly making use of Lemma \ref{asymptotic H_1
G_m}.
 \epf

We record one more lemma.

\bl \label{asymptotic compare}
 Suppose that $M$ and $N$ are two finitely generated $\Op\ps{G}$-modules
with $M = M(\pi)$ and $N= N(\pi)$. Assume that there is a
$\Op\ps{G}$-homomorphism $\varphi :M\lra N$ which has pseudo-null
kernel and cokernel. Then for $i\geq 0$, we have
\[ \ord_q \big(H_i(G_m, M)\big) = \ord_q \big(H_i(G_m, N)\big) +
O(p^{(r-1)m}) \]\el

\bpf
 The statement will follow if it holds in the two special cases of exact
sequences
\[\ba{c} 0\lra M \lra N \lra P\lra 0, \\
  0\lra P \lra M \lra N\lra 0, \ea
\] where $P$ is a pseudo-null $\Op\ps{G}$-module. Note that $P=P(\pi)$.
Taking $G_m$-homology of the first exact sequence, we have
\[ H_1(G_m, P) \lra M_{G_m} \lra N_{G_m} \lra P_{G_m} \lra 0. \]
By Lemmas \ref{asymptotic G_m} and \ref{asymptotic H_1 G_m}, we have
\[ \ord_q\big(H_1(G_m, P)\big) = \ord_q\big(P_{G_m}\big) = O(q^{(r-1)m}).
\] The cases for the higher cohomology groups and for the second exact
 sequence can be proven similarly. \epf

We can now prove our theorem.

\bpf[Proof of Theorem \ref{asymptotic formula}]
 Let
 \[\varphi: M \lra \bigoplus_{i=1}^{s}\Op\ps{G}/\pi^{\al_i}\] be
a $\Op\ps{G}$-homomorphism whose kernel and cokernel are pseudo-null
(recall that we are assuming $M= M(\pi)$). By Lemma \ref{asymptotic
compare}, we have
\[ \ord_q\big(H_i(G_m, M)\big) =
\sum_{i=1}^s \ord_q\Big(H_i\big(G_m,\Op\ps{G}/\pi^{\al_i}\big)\Big)
 + O(p^{(r-1)m}). \]
 Therefore, we are reduced to showing that
 \[ \ord_q\Big(H_i\big(G_m,\Op\ps{G}/\pi^{\al}\big)\Big)=
\begin{cases} \al p^{rm} & \text{\mbox{if}
$i =0$},\\
 0 & \text{\mbox{if} $i\geq 1$.}
\end{cases} \]
Since $\Op\ps{G}$ has no zero divisors, we have an exact sequence
 \[ 0\lra \Op\ps{G} \stackrel{\cdot \pi^\al}{\lra} \Op\ps{G}
 \lra \Op\ps{G}/\pi^{\al} \lra 0. \]
 Since $H_i(G_m, \Op\ps{G}) = 0$ for $i\geq
1$, it follows from considering the $G_m$-homology that \[H_i(G_m,
\Op\ps{G}/\pi^{\al}) =0\] for $i\geq 1$. It remains to show the
first equality, and this is immediate from the facts that
$\Op[G/G_m]/\pi^{\al} \cong (\Op/\pi^{\al})^{p^{rm}}$ (as abelian
groups) and that $|\Op/\pi^{\al}| = q^{\al}$.
 \epf

\section{Arithmetic Preliminaries} \label{Arithmetic Preliminaries}

In this section, we introduce the Selmer groups and Selmer
complexes. Here, we fix the notation that we shall use throughout
this section. To start, let $p$ be a prime. We let $F$ be a number
field. If $p=2$, we assume further that $F$ has no real primes.
Denote $\Op$ to be the ring of integers of some finite extension $K$
of $\Qp$, and fix a local parameter $\pi$ for $\Op$. Suppose that we
are given the following datum $\big(A, \{A_v\}_{v|p},
\{A^+_v\}_{v|\R} \big)$ defined over $F$:

\begin{enumerate}
 \item[(a)] $A$ is a
cofinitely generated cofree $\Op$-module of $\Op$-corank $d$ with a
continuous, $\Op$-linear $\Gal(\bar{F}/F)$-action which is
unramified outside a finite set of primes of $F$.

 \item[(b)] For each prime $v$ of $F$ above $p$, $A_v$ is a
$\Gal(\bar{F}_v/F_v)$-submodule of $A$ which is cofree of
$\Op$-corank $d_v$.

\item[(c)] For each real prime $v$ of $F$, we write $A_v^+=
A^{\Gal(\bar{F}_v/F_v)}$  which is assumed to be cofree of
$\Op$-corank $d^+_v$.

\item[(d)] The following equality
  \begin{equation} \label{data equality}
  \sum_{v|p} (d-d_v)[F_v:\Qp] = dr_2(F) +
 \sum_{v~\mathrm{real}}(d-d^+_v)
  \end{equation}
holds. Here $r_2(F)$ denotes the number of complex primes of $F$.
\end{enumerate}

 We now consider the base change property of our datum. Let $L$ be a
finite extension of $F$. We can then obtain another datum $\big(A,
\{A_w\}_{w|p}, \{A^+_w\}_{w|\R} \big)$ over $L$ as follows: we
consider $A$ as a $\Gal(\bar{F}/L)$-module, and for each prime $w$
of $L$ above $p$, we set $A_w =A_v$, where $v$ is a prime of $F$
below $w$, and view it as a $\Gal(\bar{F}_v/L_w)$-module. Then $d_w
= d_v$. For each real prime $w$ of $L$, one sets
$A^{\Gal(\bar{L}_w/L_w)}= A^{\Gal(\bar{F}_v/F_v)}$ and writes $d^+_w
= d^+_v$, where $v$ is a real prime of $F$ below $w$. In general,
the $d_w$'s and $d_w^+$'s need not satisfy equality (\ref{data
equality}). We now record the following lemma which gives some
sufficient conditions for equality (\ref{data equality}) to hold for
the datum $\big(A, \{A_w\}_{w|p}, \{A^+_w\}_{w|\R} \big)$ over $L$.

\setcounter{thm}{0}

\bl \label{data base change} Suppose that $\big(A, \{A_v\}_{v|p},
\{A^+_v\}_{v|\R} \big)$ is a datum defined over $F$.
 Suppose further that at least one of the following statements holds.
 \begin{enumerate}
\item[$(i)$] All the archimedean primes of $F$ are unramified in
$L$.

\item[$(ii)$] $[L:F]$ is odd

\item[$(iii)$] $F$ is totally imaginary.

\item[$(iv)$] $F$ is totally real, $L$ is totally imaginary and
\[ \sum_{v~\mathrm{real}} d^+_v = d[F:Q]/2.\]
 \end{enumerate}
Then we have the equality
 \[ \sum_{w|p} (d-d_w)[L_w:\Qp] = dr_2(L) +
 \sum_{w~\mathrm{real}}(d-d^+_w).\]
 \el

\bpf Note that if either of the assertions in (ii) or (iii) holds,
then the assertion in (i) holds. Therefore, to prove the lemma in
these cases, it suffices to prove it under the assumption of (i). We
first perform the following calculation
 \[ \ba{rl}
\displaystyle \sum_{w|p} (d-d_w)[L_w:\Qp]\!\!
   &= \displaystyle\sum_{v|p}\sum_{w|v} (d-d_v)[L_w:F_v][F_v:\Qp] \\
   &= \displaystyle\sum_{v|p} (d-d_v)[F_v:\Qp] \sum_{w|v}[L_w:F_v] \\
   &= \displaystyle [L:F]\sum_{v|p} (d-d_v)[F_v:\Qp] \\
  &= \displaystyle d[L:F]r_2(F) + [L:F]\sum_{v~\mathrm{real}} (d-d^+_v). \\  \ea \]
 Now if (i) holds, then every prime of $L$ above a real prime (resp.,
complex prime) of $F$ is a real prime (resp., complex prime).
Therefore, one has $[L:F]r_2(F) = r_2(L)$ and
\[ [L:F]\sum_{v~\mathrm{real}} (d-d^+_v) = \sum_{w~\mathrm{real}} (d-d^+_w). \]
The required conclusion then follows.

Now suppose that (iv) holds. Then $r_2(F) =0$ and we have
 \[ \ba{rl}
\displaystyle \sum_{w|p} (d-d_w)[L_w:\Qp] =\displaystyle
[L:F]\sum_{v~\mathrm{real}} (d-d^+_v)\!\! &=
  \displaystyle[L:F]\sum_{v~\mathrm{real}}d - [L:F]\sum_{v~\mathrm{real}}d^+_v \\
   &= [L:F][F:\Q]d - [L:F] d[F:\Q]/2 \vspace{0.1in} \\
 &= d[L:\Q]/2 = dr_2(L).
\ea
\] \epf

We now describe briefly the arithmetic situation, where we can
obtain the above set of data from. Let $V$ be a $d$-dimensional
$K$-vector space with a continuous $\Gal(\bar{F}/F)$-action which is
unramified outside a finite set of primes. Suppose that for each
prime $v$ of $F$ above $p$, there is a $d_v$-dimensional
$K$-subspace $V_v$ of $V$ which is invariant under the action of
$\Gal(\bar{F}_v/F_v)$, and for each real prime $v$ of $F$,
$V^{\Gal(\bar{F}_v/F_v)}$ has dimension $d_v^+$. Choose a
$\Gal(\bar{F}/F)$-stable $\Op$-lattice $T$ of $V$ (Such a lattice exists by
compactness). We can obtain a data as above from $V$ by setting $A=
V/T$ and $A_v = V_v/ T\cap V_v$. Note that $A$ and $A_v$ depends on
the choice of the lattice $T$. We mention some concrete examples.

\medskip
(1) Let $B$ be an abelian variety of dimension $g$ defined over a
number field $F$. For simplicity, we assume that $F$ is totally
imaginary and that the abelian variety $B$ has semistable reduction
over $F$. We define a set of data $(A, \{A_v\})$ by first setting $A
= B[p^{\infty}]$. For each prime $v$ of $F$ above $p$, let
$\mathcal{F}_v$ be the formal group attached to the Neron model for
$B$ over the ring of integers $\Op_{F_v}$ of $F_v$, and we assume
that $\mathcal{F}_v$ is a formal group of height $g$ for all $v|p$.
For instance, this latter condition is satisfied if $B$ has good
ordinary reduction at all $v|p$. We then set $A_v =
\mathcal{F}_v(\overline{\m})[p^{\infty}]$, where $\overline{\m}$ is
the maximal ideal of the rings of integers of $\overline{F}_v$. Note
that $A_v\cong (\Qp/\Zp)^g$ as an abelian group by our height
assumption. It is easy to see that $(A, \{A_v\})$ satisfies the
condition in Section \ref{Arithmetic Preliminaries} by taking $d =
2g$ and $d_v = g$. It is worthwhile mentioning that the (strict)
Selmer groups attached to this set of data coincide with the
classical Selmer groups of the abelian variety, when the Selmer
groups are considered over an admissible $p$-adic Lie extension (see
\cite{CG}).

\medskip
(2) More generally, a source of examples where we can obtain such a
datum is that of a nearly ordinary Galois representation in the
sense of Weston \cite{We}. This is a finite-dimensional $K$-vector
space equipped with a $K$-linear action of the absolute Galois group
$\Gal(\bar{F}/F)$ such that for each prime $v$ of $F$ dividing $p$,
there is sequence
\[ 0 = V_{v,0} \subseteq V_{v,1}\subseteq\cdots\subseteq V_{v,n} = V\]
 of nearly ordinary $\Gal(\bar{F}_v/F_v)$-subspace of $V$, where
$V_{v,i}$ has $K$-dimension $i$. Following \cite{We}, a set of
Selmer weights for $V$ is a choice of integers $c_v$  ($0\leq
c_v\leq d$) for each $v$ dividing $p$ such that
\[ \sum_{v|p} c_v[F_v:\Qp] = dr_2(F) +
 \sum_{v~\mathrm{real}}(d-d^+_v),
  \]
where $r_2(F)$ denotes the number of complex primes of $F$. Set $V_v
= V_{v, d-c_v}$. (In other words, our $d_v$ is $d-c_v$.) For more
concrete examples, we refer readers to \cite[\S 9]{G89} and
\cite[Section 1.2]{We}.

\subsection{Selmer groups}

We now introduce the Selmer groups. Let $S$ be a finite set of
primes of $F$ which contains all the primes above $p$, the ramified
primes of $A$ and all the infinite primes of $F$. Denote $F_S$ to be the
maximal algebraic extension of $F$ unramified outside $S$ and write
$G_S(\mathcal{L}) = \Gal(F_S/\mathcal{L})$ for every algebraic
extension $\mathcal{L}$ of $F$ which is contained in $F_S$. Let $L$
be a finite extension of $F$ contained in $F_S$ such that the data
$\big(A, \{A_w\}_{w|p}, \{A^+_w\}_{w|\R} \big)$ satisfies (\ref{data
equality}). For a prime $w$ of $L$ lying over $S$, set
\[ H^1_{str}(L_w, A)=
\begin{cases} \ker\big(H^1(L_w, A)\lra H^1(L_w, A/A_w)\big) & \text{\mbox{if} $w$
 divides $p$},\\
 \ker\big(H^1(L_w, A)\lra H^1(L^{ur}_w, A)\big) & \text{\mbox{if} $w$ does not divide $p$,}
\end{cases} \]
 where $L_w^{ur}$ is the maximal unramified extension of $L_w$.
The (strict) Selmer group attached to the data is then defined by
\[ S(A/L) := \Sel^{str}(A/L) := \ker\Big( H^1(G_S(L),A)\lra
\bigoplus_{w \in S_L}H^1_s(L_w, A)\Big),\] where we write
$H^1_s(L_w, A) = H^1(L_w, A)/H^1_{str}(L_w, A)$ and $S_L$ denotes
the set of primes of $L$ above $S$. It is straightforward to verify
that $S(A/L) = \ilim_n S(A[\pi^n]/L)$, where $S(A[\pi^n]/L)$ is the
Selmer group defined similarly as above by replacing $A$ by
$A[\pi^n]$ and $A_v$ by $A_v[\pi^n]$. Here the direct limit is taken
over the maps $S(A[\pi^n]/L) \lra S(A[\pi^{n+1}]/L)$ which are
induced by the natural injections $A[\pi^n]\hookrightarrow
A[\pi^{n+1}]$ and $A_w[\pi^n]\hookrightarrow A_w[\pi^{n+1}]$. We
will write $X(A/L)$ for its Pontryagin dual.

\smallskip
 We shall say that $F_{\infty}$ is an
\textit{$S$-admissible $p$-adic Lie extension} of $F$ if (i)
$\Gal(F_{\infty}/F)$ is compact $p$-adic Lie group, (ii)
$F_{\infty}$ contains the cyclotomic $\Zp$ extension $F_{\cyc}$ of
$F$ and (iii) $F_{\infty}$ is unramified outside $S$. Write $G =
\Gal(F_{\infty}/F)$, $H = \Gal(F_{\infty}/F_{\cyc})$ and $\Ga
=\Gal(F_{\cyc}/F)$. We define $S(A/F_{\infty}) = \ilim_L S(A/L)$,
where the limit runs over all finite extensions $L$ of $F$ contained
in $F_{\infty}$. We write $X(A/F_{\infty})$ for the Pontryagin dual
of $S(A/F_{\infty})$. By a similar argument to that in
\cite[Corollary 2.3]{CS12}, one can show that $X(A/F_{\infty})$ is
independent of the choice of $S$ as long as $S$ contains all the
primes above $p$, the ramified primes of $A$, the primes that ramify
in $F_{\infty}/F$ and all infinite primes.

We introduce another variant of the Selmer groups which is usually
called the Greenberg Selmer group. Now set
\[ H^1_{Gr}(F_v, A)=
\begin{cases} \ker\big(H^1(F_v, A)\lra H^1(F_v^{ur}, A/A_v)\big) & \text{\mbox{if} $v|p$},\\
 \ker\big(H^1(F_v, A)\lra H^1(F^{ur}_v, A)\big) & \text{\mbox{if} $v\nmid p$.}
\end{cases} \]

The Greenberg Selmer group attached to these data is then defined by
\[ \Sel^{Gr}(A/F) = \ker\Big( H^1(G_S(F),A)\lra \bigoplus_{v \in S}H^1_g(F_v,
A)\Big),\] where we write $H^1_g(F_v, A) = H^1(F_v, A)/H^1_{Gr}(F_v,
A)$. For an $S$-admissible $p$-adic Lie extension $F_{\infty}$, we
define $\Sel^{Gr}(A/F_{\infty}) = \ilim_L \Sel^{Gr}(A/L)$ and denote
$X^{Gr}(A/F_{\infty})$ to be the Pontryagin dual of
$\Sel^{Gr}(A/F_{\infty})$. The following lemma compares the two
Selmer groups of Greenberg.

\bl
  We have an exact sequence
\[0\lra S(A/F_{\infty})\lra \Sel^{Gr}(A/F_{\infty}) \lra N\lra 0,\]
where $N$ is a cofinitely generated $\Op\ps{H}$-module. \el

\bpf
 Now consider the following commutative
diagram
\[  \entrymodifiers={!! <0pt, .8ex>+} \SelectTips{eu}{}\xymatrix{
    0 \ar[r]^{} & S(A/F_{\infty}) \ar[d] \ar[r] &  H^1(G_S(F_{\infty}), A) \ar@{=}[d]
    \ar[r] & \bigoplus_{v \in S} J_v(A/F_{\infty}) \ar[d]^{\al} \\
    0 \ar[r]^{} & \Sel^{Gr}(A/F_{\infty}) \ar[r]^{} & H^1(G_S(F_{\infty}), A) \ar[r] & \
    \bigoplus_{v \in S} J^{Gr}_v(A/F_{\infty})  } \]
with exact rows, where $J^{Gr}_v(A/F_{\infty}) =
\ilim_L\bigoplus_{w|v} H^1_g(L_w, A)$. It therefore remains to show
that $\ker \al$ is cofinitely generated over $\Op\ps{H}$. Clearly,
$J_v(A/F_{\infty}) = J^{Gr}_v(A/F_{\infty})$ for $v\nmid p$. For
each $v|p$, fix a prime $w$ of $F_{\infty}$ above $v$. Write
$I_{\infty, w}$ for the inertia subgroup of
$\Gal(\overline{F}_{\infty,w}/F_{\infty, w})$ and $U_w =
\Gal(\overline{F}_{\infty,w}/F_{\infty, w})/I_{\infty, w}$. It then
follows from the Hochschild-Serre spectral sequence that we have
\[ 0\lra H^1(U_w, (A/A_v)^{I_{\infty,w}}) \lra H^1(H_w, A/A_v)
\lra H^1(I_{\infty,w}, A/A_v)^{U_w}.\]
 Since $U_w$ is topologically cyclic, $H^1(U_w,
(A/A_v)^{I_{\infty,w}}) \cong
\big((A/A_v)^{I_{\infty,w}}\big)_{U_v}$ and so is cofinitely
generated over $\Op$. Since $F_{\infty}$ contains $F^{\cyc}$, the
decomposition group of $G$ at $v$ has at least dimension one for
each $v|p$. Hence it follows that $\ker \al$ is cofinitely generated
over $\Op\ps{H}$, as required. \epf

\bl \label{Greenberg equal}
 One has  \[ \rank_{\Op\ps{G}}\big(
X(A/F_{\infty})\big) = \rank_{\Op\ps{G}}
\big(X^{Gr}(A/F_{\infty})\big),
\]
 and $X(A/F_{\infty})(\pi)$ and $X^{Gr}(A/F_{\infty})(\pi)$
have the same elementary representations.
 \el

\bpf
 By the preceding lemma, one has an exact sequence
 \[0\lra N'\lra X^{Gr}(A/F_{\infty}) \lra X(A/F_{\infty})\lra 0\]
for some finitely generated $\Op\ps{H}$-module $N'$. The first
equality is immediate, and the second assertion follows from Lemma
\ref{pseudo-isomorphism lemma}. \epf

\subsection{Selmer complexes}
\label{Selmer complexes}

We now consider the Selmer complex associated to the data $\big(A,
\{A_v\}_{v|p}, \{A^+_v\}_{v|\R}\big)$. The notion of a Selmer
complex was first conceived and introduced in \cite{Nek}. In our
discussion, we consider a modified version of the Selmer complex as
given in \cite[4.2.11]{FK}. Write $T^* =
\Hom_{\cts}(A,\mu_{p^{\infty}})$ and $T^*_v =
\Hom_{\cts}(A/A_v,\mu_{p^{\infty}})$. For every finite extension $L$
of $F$ and $w$ a prime of $L$ above $p$, write $T^*_w = T^*_v$,
where $v$ is the prime of $F$ below $w$. For any profinite group
$\mathcal{G}$ and a topological abelian group $M$ with a continuous
$\mathcal{G}$-action, we denote by $C(\mathcal{G}, M)$ the complex
of continuous cochains of $\mathcal{G}$ with coefficients in $M$.
Let $F_{\infty}$ be an $S$-admissible extension of $F$ with Galois
group $G$. We define a $(\Op\ps{G})[G_S(F)]$-module
$\mathcal{F}_G(T^*)$ as follows: as an $\Op$-module,
$\mathcal{F}_G(T^*) = \Op\ps{G}\ot_{\Op}T^*$, and the action of
$G_S(F)$ is given by the formula $\sigma(x\ot t) =
x\bar{\sigma}^{-1}\ot \sigma t$, where $\bar{\sigma}$ is the
canonical image of $\sigma$ in $G \subseteq \Op\ps{G}$. We define
the $(\Op\ps{G})[\Gal(\bar{F}_v/F_v)]$-module $\mathcal{F}_G(T_v^*)$
in a similar fashion.

For every prime $v$ of $F$, we write $C\big(F_v,
\mathcal{F}_G(T^*)\big)= C\big(\Gal(\bar{F}_v/F_v),
\mathcal{F}_G(T^*)\big)$. For each prime $v$ not dividing $p$,
denote $C_f\big(F_v, \mathcal{F}_G(T^*)\big)$ to be the subcomplex
of $C\big(F_v, \mathcal{F}_G(T^*)\big)$, whose degree $m$-component
is $0$ unless $m\neq  0, 1$, whose degree $0$-component is
$C^0\big(F_v, \mathcal{F}_G(T^*)\big)$, and whose degree
$1$-component is
\[ \ker\Big( C^1\big(F_v, \mathcal{F}_G(T^*)\big)_{d=0} \longrightarrow
H^1\big(F_v^{ur},\mathcal{F}_G(T^*)\big) \Big). \]

The Selmer complex $SC(T^*, T^*_v)$ is then defined to be
\[ \mathrm{Cone}\bigg( C\big(G_S(F), \mathcal{F}_G(T^*)\big) \longrightarrow
\displaystyle\bigoplus_{v|p}C\big(F_v,
\mathcal{F}_G(T^*)/\mathcal{F}_G(T^*_v)\big)\oplus \displaystyle
\bigoplus_{v\nmid p} C\big(F_v, \mathcal{F}_G(T^*)\big)/C_f\big(F_v,
\mathcal{F}_G(T^*)\big)\bigg)[-1].
\]
 Here $[-1]$ is the translation by $-1$ of the complex. We will
write  $H^i\big(SC(T^*, T^*_v)\big)$ for the $i$th cohomology group
of the complex $SC(T^*, T^*_v)$. We now state the following
proposition which is proven in \cite[Proposition 4.2.35]{FK}.

\bp \label{FK prop} Let $\mathcal{G}$ be the kernel of
$\Gal(\bar{F}/F) \lra G$. For a place $v$ of $F$, fixing an
embedding $F \hookrightarrow F_v$, let $\mathcal{G}(v)$ be the
kernel of $\Gal(\bar{F}_v/F_v)\lra G$ and let $G_v\subseteq G$ be
the image. Then the following statements hold.
\begin{enumerate}
\item[$(a)$] $H^i\big(SC(T^*, T^*_v)\big) = 0$ for $i\neq 1, 2,
3$.

\item[$(b)$] We have an exact sequence
\[ \ba{c}
0 \lra X(A/F_{\infty}) \lra H^2\big(SC(T^*, T^*_v)\big) \lra
\displaystyle \bigoplus_{v|p} \Op\ps{G}\ot_{\Op\ps{G_v}}
\big(T^*_v(-1)\big)_{\mathcal{G}(v)}
 \hspace{2in}\\
 \hspace{1.5 in} \lra \big(T^*(-1)\big)_{\mathcal{G}} \lra H^3\big(SC(T^*,
 T^*_v)\big) \lra 0.
 \ea\]
\end{enumerate}
\ep

Since $F_{\infty}$ contains $F_{\cyc}$, it follows that for every
prime $v|p$, the group $G_v$ has dimension at least 1.  Therefore,
$\bigoplus_{v|p} \Op\ps{G}\ot_{\Op\ps{G_v}}
\big(T^*_v(-1)\big)_{\mathcal{G}(v)}$ is finitely generated over
$\Op\ps{H}$, and one can apply Lemma \ref{pseudo-isomorphism lemma}
to obtain the following statement.

\bl \label{Greenberg equal complex}
 $X(A/F_{\infty})(\pi)$ and $H^2\big(SC(T^*, T^*_v)\big)(\pi)$
have the same elementary representations.
 \el

We end the section with the following remark.

\br It is clear from the exact sequence in Proposition \ref{FK prop}
that $H^3\big(SC(T^*, T^*_v)\big)$ is finitely generated over $\Op$.
In particular, this implies that $H^3\big(SC(T^*, T^*_v)\big)(\pi)$
is pseudo-null. Now assume further that $X(A/F_{\infty})$ satisfies
the $\M_H(G)$-property in the sense of \cite{CFKSV, CS12, FK} (see
also Subsection \ref{MHG subsection}). Note that when $F_{\infty}$
is of dimension 1, this is equivalent to saying that
$X(A/F_{\infty})$ is a torsion module. Then it follows from
\cite[Proposition 4.3.11]{FK} that $H^1\big(SC(T^*, T^*_v) = 0$
(resp., $H^1\big(SC(T^*, T^*_v)$ is finitely generated over $\Op$)
if $F_{\infty}$ is of dimension $>1$ (resp., dimension 1). In either
cases, we have that $H^1\big(SC(T^*, T^*_v)\big)(\pi)$ is
pseudo-null. Therefore, it follows from the above discussion,
Proposition \ref{FK prop}(a) and Lemma \ref{Greenberg equal complex}
that the $\pi$-primary submodule of $X(A/F_{\infty})$ essentially
captures the ``$\pi$-primary component" of the Selmer complex
$SC(T^*, T^*_v)$. \er

\section{$\pi$-submodules of dual Selmer groups}
\label{pi-submodules of dual Selmer groups}

Throughout this section, $\big(A, \{A_v\}_{v|p}, \{A^+_v\}_{v|\R}
\big)$ is a datum defined as in Section 3 over a number field $F$.
As before, $S$ will denote a finite set of primes that contain the
primes of $F$ above $p$, the ramified primes of $A$ and the
archimedean primes of $F$. Let $F_{\infty}$ be an $S$-admissible
$p$-adic Lie extension of $F$ whose Galois group
$G=\Gal(F_{\infty}/F)$ is a pro-$p$ torsion-free $p$-adic Lie group of dimension $r$.
We also recall that $G_0$ is a fixed open normal uniform subgroup of $G$, and $G_m$ is denoted to be the $(m+1)$-term of the
lower $p$-series of $G_0$ which is defined by
\[ G_{m+1} = G_m^{p}[G_m,G], ~\mbox{for}~ m\geq 0. \]
 Denote $F_m$ to be the fixed field of $G_m$. Note that this is a finite Galois
extension of $F$ of degree $[G:G_0]p^{rm}$.

For the remainder of the paper, we will work with $X(A/F_{\infty})$.
In view of Lemmas \ref{Greenberg equal} and \ref{Greenberg equal
complex}, all our main results (Theorem \ref{tate dual} and Theorem
\ref{congruent}) also hold for the Greenberg Selmer groups and the
second cohomology of the Selmer complexes as defined in Subsection
\ref{Selmer complexes}.

If $G$ is a pro-$p$ group, we write $h_1(G) = \dim_{\Z/p}\big(H^1(G,
\Z/p)\big)$ and $h_2(G) = \dim_{\Z/p}\big(H^2(G, \Z/p)\big)$. If $M$
is a cofinitely generated $\Op$-module, we denote $M_{\mathrm{div}}$
to be the maximal $\Op$-divisible submodule of $M$. We now record a
useful lemma which allows us to estimate the order of certain
cohomology groups, whose easy proof is left to the reader.

\setcounter{thm}{0}

\bl \label{cohomology order inequalities} Let $G$ be a pro-$p$
group, and let $M$ be a discrete $G$-module which is cofinitely
generated over $\Op$. Let $n$ be a positive integer.
 If $h_1(G)$ is finite, then for every $n\geq 1$, $H^1(G,M)[\pi^n]$ is finite
 and
  \[
 \ord_q\big(H^1(G,M)[\pi^n]\big) \leq n h_1(G) \big(\mathrm{corank}_{\Op}(M) +
\ord_q( M/M_{\mathrm{div}})\big).
\]

If $h_2(G)$ is finite, then for every $n\geq 1$, $H^2(G,M)[\pi^n]$
is finite and
 we have the following inequality
 \[
 \ord_q\big(H^2(G,M)[\pi^n]\big) \leq n h_2(G)\big(\mathrm{corank}_{\Op}(M) +
\ord_q( M/M_{\mathrm{div}})\big).
\] \el

\subsection{Tate dual}

In this subsection, we prove our first main result of the paper. For
a given set of data $\big(A, \{A_v\}_{v|p}, \{A^+_v\}_{v|\R} \big)$,
we define its (Tate) dual data as follows. For a $\Op$-module $N$,
we denote $T_{\pi}(N)$ to be its $\pi$-adic Tate module, i.e.,
$T_{\pi}(N) = \plim_i N[\pi^i]$.  We then set $A^* =
\Hom_{\cts}(T_{\pi}(A),\mu_{p^{\infty}})$. Similarly, for each $v|p$
(resp., $v$ real), we set $A^*_v=
\Hom_{\cts}(T_{\pi}(A/A_v),\mu_{p^{\infty}})$ (resp., $(A^*)^+_v=
\Hom_{\cts}(T_{\pi}(A/A^+_v),\mu_{p^{\infty}})$).  It is an easy
exercise to verify that $\big(A^*, \{A_v^*\}_{v|p},
\{(A^*)^+_v\}_{v|\R} \big)$ satisfies equality (\ref{data
equality}). Therefore, we can attach Selmer groups to this dual data
which we denote by $S(A^*/F_{\infty})$ and
$S(A^*[\pi^n]/F_{\infty})$. We then denote $X(A^*/F_{\infty})$ to be
the Pontryagin dual of $S(A^*/F_{\infty})$. We are now in the
position to state the first main theorem of the paper.

\bt \label{tate dual} Let $F_{\infty}$ be an admissible $p$-adic Lie
extension such that $G=\Gal(F_{\infty}/F)$ is pro-$p$ torsion-free $p$-adic Lie group. Then $X(A/F_{\infty})$ and $X(A^*/F_{\infty})$ have the same
$\Op\ps{G}$-ranks, and $X(A/F_{\infty})(\pi)$ and
$X(A^*/F_{\infty})(\pi)$ have the same elementary representations.
\et

For data coming from (nearly) ordinary representations, it is
expected that $X(A/F_{\infty})$ is a torsion $\Op\ps{G}$-module (see
\cite[Conjecture 1]{G89} or \cite[Conjecture 1.7]{We}). We therefore
record the following important corollary.

\bc \label{tate dual corollary}
 Let $F_{\infty}$ be an admissible $p$-adic Lie
extension such that $G=\Gal(F_{\infty}/F)$ is a pro-$p$ torsion-free $p$-adic Lie group. Then $X(A/F_{\infty})$ is a torsion $\Op\ps{G}$-module if and
only if $X(A^*/F_{\infty})$ is a torsion $\Op\ps{G}$-module. \ec

The remainder of the subsection will be devoted to the proof of
Theorem \ref{tate dual}. By Proposition \ref{pseudo-isomorphic3}, we
are reduced to proving the following proposition.

\bp \label{mu tate dual}
 For every $n\geq 1$, we have
 \[
 \mu_{\Op\ps{G}}\Big(X(A/F_{\infty})/\pi^n\Big)=
\mu_{\Op\ps{G}}\Big(X(A^*/F_{\infty})/\pi^n\Big).
\]
\ep

\noindent \textit{Proof.}  Let $n$ be an arbitrary fixed positive
integer. Then for every $m\geq 1$, we have
 \[ \ba{l}
  \bigg| [G:G_0]\Big(\mu_{\Op\ps{G}}\big(X(A/F_{\infty})/\pi^n\big)-
\mu_{\Op\ps{G}}\big(X(A^*/F_{\infty})/\pi^n\big)\Big)p^{rm} \bigg|
\leq \\
 \hspace{1in} \Big| [G:G_0]\mu_{\Op\ps{G}}\big(X(A/F_{\infty})/\pi^n\big)p^{rm} -
 \ord_q\big(S(A/F_{\infty})[\pi^n]^{G_m}\big)\Big| \\
 \hspace{1in} +
  \Big| [G:G_0]\mu_{\Op\ps{G}}\big(X(A^*/F_{\infty})/\pi^n\big)p^{rm} -
 \ord_q\big(S(A^*/F_{\infty})[\pi^n]^{G_m}\big)\Big| \\
 \hspace{1in} + \Big| \ord_q \big(S(A[\pi^n]/F_m)\big) -
 \ord_q\big(S(A/F_{\infty})[\pi^n]^{G_m}\big)\Big| \\
 \hspace{1in} +
 \Big| \ord_q \big(S(A^*[\pi^n]/F_m)\big) -
 \ord_q\big(S(A^*/F_{\infty})[\pi^n]^{G_m}\big)\Big| \\
 \hspace{1in} + \Big| \ord_q \big(S(A[\pi^n]/F_m)\big) -
 \ord_q\big(S(A^*[\pi^n]/F_m)\big)\Big|.
 \ea
\]
 The required equality of the proposition will follow once we can show that each of the five
quantities on the right is $O(p^{(r-1)m})$. The first two quantities
are $O(p^{(r-1)m})$ by Theorem \ref{asymptotic formula}. We now
proceed to show that the third and fourth quantities are
$O(p^{(r-1)m})$. To show this, we first need to estimate the order
of the kernels and cokernels of the maps
 \[ S(A[\pi^n]/F_m) \stackrel{r_m}{\lra} S(A/F_m)[\pi^n]
 \stackrel{s_m}{\lra}  \big(S(A/F_{\infty})[\pi^n]\big)^{G_m}.\]
 One sees easily that $\ker r_m \subseteq A(F_m)/\pi^n$ and $\ker s_m
 \subseteq H^1(G_m, A(F_{\infty}))[\pi^n]$. It is clear that one has
 $\ord_q(\ker r_m) \leq nr$ for every $m$, and therefore,
 $\ord_q(\ker r_m) = O(1)$. On the other hand, it follows from Lemma \ref{cohomology order inequalities}
that $\ord_q\big(H^1(G_m, A(F_{\infty}))[\pi^n]\big) = O(1)$ (noting
that $h_1(G_m)$ is a constant function in $m$). Thus, one has
$\ord_q(\ker s_m) =  O(1)$.

To estimate $\coker r_m$ and $\coker s_m$, one first observes that
$\ord_q(\coker r_m) \leq \ord_q(\ker r_m')$ and that $\ord_q(\coker
s_m) \leq \ord_q(\ker s_m') + \ord_q(H^2(G_m,
A(F_{\infty}))[\pi^n])$, where $r_m'$ and $s_m'$ are given by
\[r_m' = \big(r'_{m,v_m}\big): \bigoplus_{v_m\in S_{F_m}}H^1_s(F_{m,v_m}, A[\pi^n])
\lra \bigoplus_{v_m\in S_{F_m}}H^1_s(F_{m,v_m}, A)[\pi^n];\]
\[s_m'  = \big(s'_{m,v_m}\big) : \bigoplus_{v_m\in S_{F_m}}H^1_s(F_{m,v_m}, A)[\pi^n]
\lra \Big(\ilim_m\bigoplus_{v_m\in S_{F_m}}H^1_s(F_{m,v_m},
A)[\pi^n]\Big)^{G_m}.\]
 By Lemma \ref{cohomology order inequalities}, one has that
$\ord_q\big(H^2(G_m, A(F_{\infty}))[\pi^n]\big) = O(1)$ (noting that
$h_2(G_m)$ is a constant function in $m$; in fact, one has $h_2(G_m)
= r(r-1)/2$ by \cite[Theorem 4.35]{DSMS}). To estimate $\coker
r'_m$, we first observe that

\[\ker r'_{m,v_m} \subseteq
\begin{cases} \ker\Big(H^1(F_{m, v_m}, A/A_{v_m}[\pi^n])\lra
 H^1(F_{m, v_m}, A/A_{v_m})[\pi^n]\Big)
 & \text{\mbox{if} $v_m|p$},\\
\ker \Big(H^1(F^{ur}_{m, v_m}, A[\pi^n])\lra
 H^1(F^{ur}_{m, v_m}, A)[\pi^n]\Big)  & \text{\mbox{if} $v_m\nmid p$,}
\end{cases} \]

\[  =
\begin{cases} A/A_{v_m}(F_{m, v_m})[\pi^n] \hspace{0.6in}
 & \text{\mbox{if} $v_m|p$}, \hspace{1in}\\
 A(F^{ur}_{m, v_m})[\pi^n] \hspace{0.6in}
 & \text{\mbox{if} $v_m\nmid p$, \hspace{1in}}
\end{cases} \]

It is now clear from the above that $\ord_p(\ker r'_{m,v_m})$ is
bounded independently of $m$ and $v_m$ (for a fixed $n$). Combining
these estimates with the fact that the decomposition group of $v$ in
$G$ has dimension $\geq 1$ for every $v\in S$ (since $F_{\infty}$
contains $F^{\cyc}$), one then has the estimate $\ord_q \big(\ker
r'_m\big) = O(p^{(r-1)m})$.

 To estimate $\coker
s'_m$, we observe that

\[\ker s'_{m,v_m} \subseteq
\begin{cases} H^1\big(G_{m, v_m}, A/A_{v_m}(F_{m,v_m})\big)[\pi^n]
 & \text{\mbox{if} $v_m|p$},\\
H^1\big(\Gal(F_{\infty, v_m}/F_{m,v_m}^{ur}),
A(F_{m,v_m})\big)[\pi^n] & \text{\mbox{if} $v_m\nmid p$.}
\end{cases} \]

By appealing to Lemma \ref{cohomology order inequalities}, one
verifies easily that $\ker s'_{m,v_m}$ is bounded independent of $m$
and $v_m$ (for a fixed $n$). As before, combining these estimates
with the fact that the decomposition group of $v$ in $G$ has
dimension $\geq 1$ for every $v\in S$ (since $F_{\infty}$ contains
$F^{\cyc}$), we obtain
  $\ord_q \big(\ker s'_m\big)
= O(p^{(r-1)m})$. In conclusion, we have
 \[ \ord_q \big(S(A[\pi^n]/F_m)\big) =
 \ord_q\big(S(A/F_{\infty})[\pi^n]^{G_m}\big) +O(p^{(r-1)m}).
  \]
By a similar argument, one also has
 \[ \ord_q \big(S(A^*[\pi^n]/F_m)\big) =
 \ord_q\big(S(A^*/F_{\infty})[\pi^n]^{G_m}\big) +O(p^{(r-1)m}).
  \]
Therefore, we have shown that the third and fourth quantities are
$O(p^{(r-1)m})$. For the estimate of the final quantity, we require
the following lemma.

\bl \label{ratio of Selmer} For every $n$ and $m$, we have
 \[ \frac{|S(A[\pi^n]/F_m)|}{|H^0(G_S(F_m),A[\pi^n])|}\times
 \prod_{v_m|p}|H^0(F_{m, v_m},A[p^n]/A_{v_m}[\pi^n])| \hspace{2in}\]
 \[ \hspace{2in} = \frac{|S(A^*[\pi^n]/F_m)|}
 { |H^0(G_S(F_m),A^*[\pi^n])|} \times
 \prod_{v_m|p}|H^0(F_{m,v_m}, A^*[\pi^n]/A_{v_m}^*[\pi^n])|,  \]
 where the product is taken over all the primes of $F_m$ above $p$.
\el

\bpf (Sketch of the proof)
 This is proven in the same way as \cite[Formula (53)]{G89} and
we give the general idea behind the calculations, leaving the
details to the readers. By appealing to the fact that our datum
satisfies equality (\ref{data equality}), one can verify the
following
 \[ \frac{|H^0(G_S(F_m),A[\pi^n])|
 |H^2(G_S(F_m),A[\pi^n])|}{|H^1(G_S(F_m),A[\pi^n])|}
 = \prod_{v_m|p} \frac{|H^0(F_{m,v_m},A/A_{v_m}[\pi^n])|
 |H^2(F_{m,v_m}, A/A_{v_m}[\pi^n])|}{|H^1(F_{m, v_m}v,A/A_{v_m}[\pi^n])|}  \]
  by a global-local Euler characteristic argument.
The required equality of the lemma will follow by combining the
above with a Poitou-Tate duality argument.  \epf

We continue the proof of our main theorem.

\bpf[Proof of Theorem \ref{tate dual} $($cont'd$)$]
 Clearly, the quantities $|H^0(G_S(F_m),A[p^n])|$, $|H^0(G_S(F_m),A^*[p^n])|$,
 \linebreak
$|H^0(F_{m,v_m}, A[\pi^n]/A_v[\pi^n])|$ and $|H^0(F_{m,v_m},
A^*[\pi^n]/A_v^*[\pi^n])|$ are bounded independently of $m$ and
$v_m$ (for a fixed $n$). Since there are only finite number of
primes of $F^{\cyc}$ above $p$, the decomposition group of $v$ in
$G$ has at least dimension $1$. Therefore, it follows that
$\prod_{v_m|p}|H^0(F_{m,v_m}, A[\pi^n]/A_v[\pi^n])|$ and
$\prod_{v_m|p}|H^0(F_{m,v_m}, A^*[\pi^n]/A_v^*[\pi^n])|$ are both
$q^{O(p^{(r-1)m})}$. Therefore, in conclusion, we have
 \[ \ord_q \big(S(A[\pi^n]/F_m)\big) =
 \ord_q\big(S(A^*[\pi^n]/F_m)\big)
+O(p^{(r-1)m}), \]
 as required. This completes the proof of the theorem.
\epf

\br (1) If $F_{\infty}$ is a general $p$-adic Lie extension of $F$
(that does not contain $F^{\cyc}$) which has the property such that
for each prime $v\in S$, the decomposition group of
$\Gal(F_{\infty}/F)$ at $v$ has dimension $\geq 1$, then the
argument of Theorem \ref{tate dual} carries over to give the same
conclusion.

(2) When $F_{\infty} = F^{\cyc}$, Greenberg claimed that
$X^{Gr}(A/F^{\cyc})$ and $X^{Gr}(A^*/F^{\cyc})$ might be
pseudo-isomorphic up to an $\iota$-twist (see \cite[P. 130, Equation
(66)]{G89}) and gave some examples where this pseudo-isomorphism can
be shown (see discussion after \cite[P. 130, Equation (66)]{G89}).
In view of Lemma \ref{Greenberg equal}, Theorem \ref{tate dual} may
therefore be viewed as providing a positive answer to the
$\pi$-primary part of the assertion of Greenberg. (In fact, our
result also establishes higher analog of this.) \er

\subsection{Congruent Galois representations} \label{Congruent Galois
representations}

As before, let $F_{\infty}$ be an admissible $p$-adic Lie-extension of $F$
whose Galois group is a pro-$p$ torsion-free $p$-adic Lie group of dimension $r$. We
write $G = \Gal(F_{\infty}/F)$. To state our result, we introduce
another datum $\big(B, \{B_v\}_{v|p}, \{B^+_v\}_{v|\R}\big)$ which
satisfies the conditions (a)--(d) as in Section 3. To compare the
Selmer groups, we need to expand the set $S$ of primes to contain
the ramified primes of $B$. We introduce the following important
congruence condition on $A$ and $B$ which allows us to be able to
compare the Selmer groups of $A$ and $B$.

\medskip \noindent
$\mathbf{(Cong_n)}$ : There is an isomorphism $A[\pi^{n}]\cong
B[\pi^{n}]$ of $G_S(F)$-modules which induces a
$\Gal(\bar{F}_v/F_v)$-isomorphism $A_v[\pi^{n}]\cong B_v[\pi^{n}]$
for every $v|p$.

\medskip

Clearly, $\mathbf{(Cong_n)}$ implies $\mathbf{(Cong_i)}$ for $i\leq
n$. To simplify notation, we will write $\theta_G(A) =
\theta_G\big(X(A/F_{\infty})\big)$ and $\theta_G(B) =
\theta_G\big(X(B/F_{\infty})\big)$. The following is the second main
theorem of the paper.

\bt \label{congruent} Let $F_{\infty}$ be an admissible $p$-adic Lie
extension of $F$ whose Galois group is a pro-$p$ torsion-free $p$-adic Lie group . Suppose that $\mathbf{(Cong_{\theta_G(A)+1})}$ holds
and suppose that $X(A/F_{\infty})$ is torsion over $\Op\ps{G}$. Then
$X(B/F_{\infty})$ is torsion over $\Op\ps{G}$, and
$X(A/F_{\infty})(\pi)$ and $ X(B/F_{\infty})(\pi)$ have the same
elementary representations. \et

\bpf
 By Proposition \ref{pseudo-isomorphic2}, it suffices to show that
 \[ \mu_{\Op\ps{G}}\Big(X(A/F_{\infty})/\pi^{n}\Big) =
\mu_{\Op\ps{G}}\Big(X(B/F_{\infty})/\pi^{n}\Big)\] for $1\leq n \leq
\theta_G(A)+1$. Fix such an arbitrary $n$. Then for $m\geq 1$, we
have
\[ \ba{l}
  \bigg| [G:G_0]\Big(\mu_{\Op\ps{G}}\big(X(A/F_{\infty})/\pi^n\big)-
\mu_{\Op\ps{G}}\big(X(B/F_{\infty})/\pi^n\big)\Big)p^{rm} \bigg|
\leq \\
 \hspace{1in} \Big| [G:G_0]\mu_{\Op\ps{G}}\big(X(A/F_{\infty})/\pi^n\big)p^{rm} -
 \ord_q\big(S(A/F_{\infty})[\pi^n]^{G_m}\big)\Big| \\
 \hspace{1in} +
  \Big| [G:G_0]\mu_{\Op\ps{G}}\big(X(B/F_{\infty})/\pi^n\big)p^{rm} -
 \ord_q\big(S(B/F_{\infty})[\pi^n]^{G_m}\big)\Big| \\
 \hspace{1in} + \Big| \ord_q \big(S(A[\pi^n]/F_m)\big) -
 \ord_q\big(S(A/F_{\infty})[\pi^n]^{G_m}\big)\Big| \\
 \hspace{1in} +
 \Big| \ord_q \big(S(B[\pi^n]/F_m)\big) -
 \ord_q\big(S(B/F_{\infty})[\pi^n]^{G_m}\big)\Big| \\
 \hspace{1in} + \Big| \ord_q \big(S(A[\pi^n]/F_m)\big) -
 \ord_q\big(S(B[\pi^n]/F_m)\big)\Big|.
 \ea
\]
 As seen from the argument in the proof of Theorem \ref{tate dual},
 the first four quantities on the right of the inequality
 are $O(p^{(r-1)m})$. It remains to estimate the last quantity. By the discussion
 before this theorem, we have that $\mathbf{(Cong_n)}$ holds for $1\leq n \leq
\theta_G(A)+1$. This in turn implies that
\[ S(A[\pi^{n}]/F_m) \cong S(B[\pi^{n}]/F_m)\]
for all $m$. In particular, the final quantity on the right of the
inequality is zero. Hence we have that
\[  \bigg| [G:G_0]\Big(\mu_{\Op\ps{G}}\big(X(A/F_{\infty})/\pi^n\big)-
\mu_{\Op\ps{G}}\big(X(B/F_{\infty})/\pi^n\big)\Big)p^{rm} \bigg| =
O(p^{(r-1)m})\] which implies that
\[ \mu_{\Op\ps{G}}\Big(X(A/F_{\infty})/\pi^{n}\Big) =
\mu_{\Op\ps{G}}\Big(X(B/F_{\infty})/\pi^{n}\Big),\]
 as required.  \epf

\br If $F_{\infty}$ is a general $p$-adic Lie extension of $F$ (that
does not contain $F^{\cyc}$) which has the property such that for
each prime $v\in S$, the decomposition group of $\Gal(F_{\infty}/F)$
at $v$ has dimension $\geq 1$, then the argument of Theorem
\ref{congruent} carries over to give the same conclusion. \er

\section{Miscellaneous} \label{Miscellaneous}

\subsection{Some remarks on the $\M_H(G)$-property} \label{MHG subsection}

Let $F_{\infty}$ be an admissible $p$-adic Lie extension. As before,
we write $G= \Gal(F_{\infty}/F)$, $H= \Gal(F_{\infty}/F^{\cyc})$ and
$\Ga=\Gal(F^{\cyc}/F)$. We say that an $\Op\ps{G}$-module $M$
\textit{satisfies the $\M_H(G)$-property} if $M_f:=M/M(\pi)$ is
finitely generated over $\Op\ps{H}$. It has been conjectured for
certain Galois representations coming from abelian varieties with
good ordinary reduction at $p$ or cuspidal eigenforms with good
ordinary reduction at $p$, the dual Selmer group associated to such
a Galois representation satisfies the $\M_H(G)$-property (see
\cite{CFKSV, CS12, FK}).

For the remainder of this subsection, we will assume that $G$ is a
pro-$p$ group of dimension 2 and has no elements of order $p$. As
before, $\big(A, \{A_v\}_{v|p}, \{A^+_v\}_{v|\R} \big)$ denotes a
set of data as defined in Section 3. In preparation for further
discussion, we record the following lemma which has a similar proof
to that in \cite[Corollary 3.2]{CS12}.

\bl \label{MHG lemma}
 Let $F_{\infty}$ be an $S$-admissible
$p$-adic Lie extension whose Galois group is a pro-$p$ group of
dimension 2 and has no elements of order $p$. Suppose that
$A(F^{\cyc})$ is finite. Then the following statements are
equivalent.

 \begin{enumerate}
 \item[$(a)$]  $X(A/F_{\infty})$ satisfies the $\M_H(G)$-property.
 \item[$(b)$] $X(A/F^{\cyc})$ is a torsion $\Op\ps{\Ga}$-module,
 $X(A/F_{\infty})$ is a  torsion $\Op\ps{G}$-module and
 \[ \mu_{\Op\ps{G}}\big(X(A/F_{\infty})\big) =
\mu_{\Op\ps{\Ga}}\big(X(A/F^{\cyc})\big).\] \end{enumerate} \el

We should mention that the finiteness condition on $A(F^{\cyc})$ has
been verified in many cases, and therefore, the discussion in this
subsection may apply to these situations. In the case of an abelian
variety with good ordinary reduction at $p$, this is verified in
\cite{Im} and for a cuspidal eigenform with good ordinary reduction
at $p$, this is done in \cite[Proof of Lemma 2.2]{Su}. For a more
general result on the finiteness condition for $A$ arising from the
Galois representation attached to an \'etale $i$th-cohomology group
(for $i$ odd) of a smooth proper variety with potentially good
reduction, we refer readers to \cite{CSW, KT}.

We now state the next result which compares the structural
properties of $X(A/F_{\infty})$ and $X(A^*/F_{\infty})$, where
$X(A/F_{\infty})$ is the Selmer group associated to the set of data
$\big(A, \{A_v\}_{v|p}, \{A^+_v\}_{v|\R} \big)$ and
$X(A^*/F_{\infty})$ is the Selmer group associated to $\big(A^*,
\{A^*_v\}_{v|p}, \{(A^*)^+_v\}_{v|\R} \big)$.

\bp Let $F_{\infty}$ be an admissible $p$-adic Lie extension of $F$,
whose Galois group is a pro-$p$ group of dimension 2 and has no
elements of order $p$. Furthermore, suppose that $A(F^{\cyc})$ and
$A^*(F^{\cyc})$ are finite. Then $X(A/F_{\infty})$ satisfies the
$\M_H(G)$-property if and only if $X(A^*/F_{\infty})$ satisfies the
$\M_H(G)$-property. \ep

\bpf
 It suffices to show that if $X(A/F_{\infty})$ satisfies the
$\M_H(G)$-property, then $X(A^*/F_{\infty})$ also satisfies the
$\M_H(G)$-property. Suppose that $X(A/F_{\infty})$ satisfies the
$\M_H(G)$-property. Then by Lemma \ref{MHG lemma}, we have that
$X(A/F^{\cyc})$ is a torsion $\Op\ps{\Ga}$-module,
 $X(A/F_{\infty})$ is a  torsion $\Op\ps{G}$-module and
 \[ \mu_{\Op\ps{G}}\big(X(A/F_{\infty})\big) =
\mu_{\Op\ps{\Ga}}\big(X(A/F^{\cyc})\big).\]
 By virtue of Theorem \ref{tate
dual}, we then have that $X(A^*/F^{\cyc})$ is a torsion
$\Op\ps{\Ga}$-module,
 $X(A^*/F_{\infty})$ is a  torsion $\Op\ps{G}$-module and
 \[ \mu_{\Op\ps{G}}\big(X(A^*/F_{\infty})\big) =
\mu_{\Op\ps{\Ga}}\big(X(A^*/F^{\cyc})\big).\]
 By appealing to Lemma \ref{MHG lemma} again, this in turn implies that
$X(A^*/F_{\infty})$ satisfies the $\M_H(G)$-property. \epf

We also have a similar result as above for congruent
representations. Let $\big(B, \{B_v\}_{v|p}, \{B^+_v\}_{v|\R} \big)$
be another set of data defined as in Section 3. The next proposition
can be proven similarly by combining Theorem \ref{congruent} and
Lemma \ref{MHG lemma}.

\bp \label{MHG congruent} Let $F_{\infty}$ be an admissible $p$-adic
Lie extension of $F$, whose Galois group is a pro-$p$ group of
dimension 2 and has no elements of order $p$. Assume that
$A(F^{\cyc})$ and $B(F^{\cyc})$ are finite.  Suppose that
$\mathbf{(Cong_{\theta+1})}$ holds, where $\theta =
\max\{\theta_{\Op\ps{\Ga}}(A), \theta_{\Op\ps{G}}(A)\}$. Then if
$X(A/F_{\infty})$ satisfies the $\M_H(G)$-property, so does
$X(B/F_{\infty})$. \ep

\subsection{Comparing specializations of a big Galois representation}
\label{big galois compare section}

We apply the main result in Subsection \ref{Congruent Galois
representations} to compare the Selmer groups of specializations of
a big Galois representation. As before, let $p$ be a prime. We let
$F$ be a number field. If $p=2$, we assume further that $F$ has no
real primes. Denote $\Op$ to be the ring of integers of some finite
extension $K$ of $\Qp$. We write $R=\Op\ps{T}$ for the power series
ring in one variable. Suppose that we are given the following set of
data:

\begin{enumerate}
 \item[(a)] $\A$ is a
cofinitely generated cofree $R$-module of $R$-corank $d$ with a
continuous, $R$-linear $\Gal(\bar{F}/F)$-action which is unramified
outside a finite set of primes of $F$.

 \item[(b)] For each prime $v$ of $F$ above $p$, $\A_v$ is a
$\Gal(\bar{F}_v/F_v)$-submodule of $\A$ which is cofree of
$R$-corank $d_v$.

\item[(c)] For each real prime $v$ of $F$, we write $\A_v^+=
\A^{\Gal(\bar{F}_v/F_v)}$ which we assume to be cofree of $R$-corank
$d_v^+$.

\item[(d)] The following equality
\[
  \sum_{v|p} (d-d_v)[F_v:\Qp] = dr_2(F) +
 \sum_{v~\mathrm{real}}(d-d^+_v)
  \]
holds. Here $r_2(F)$ denotes the number of complex primes of $F$.
\end{enumerate}

For any prime element $f$ of $\Op\ps{T}$ such that $\Op\ps{T}/f$ is
a maximal order, then we can obtain a data $\big(\A[f],
\{\A_v[f]\}_{v|p}, \{\A^+_v[f]\}_{v|\R}\big)$ in the sense of
Section \ref{Arithmetic Preliminaries}. The next lemma has a easy
proof which is left to reader.

\bl \label{big galois lemma}
 Let $f$ and $g$ be prime elements of $\Op\ps{T}$ with $\pi^n|(f-g)$
 such that $\Op\ps{T}/f$ and $\Op\ps{T}/g$ are maximal orders.
 Then $\A[f, \pi^n] = \A[g,\pi^n]$. One also has similar conclusions
 for $\A_v$ and $\A^+_v$.
\el

The next proposition compares the $\pi$-primary submodules of the
dual Selmer groups of various specializations of a big Galois
representation. For a real number $x$, we denote $\lceil x \rceil$
to be the smallest integer not less than $x$.

\bp \label{big galois compare}
 Let $F_{\infty}$ be an admissible $p$-adic Lie extension of $F$
such that $G =\Gal(F_{\infty}/F)$ is uniform pro-$p$ group. Let $f$
be a prime element of $\Op\ps{T}$ such that $\Op':=\Op\ps{T}/f$ is a
maximal order. Set $A = \A[f]$ and suppose that $X(A/F_{\infty})$ is
torsion over $\Op'\ps{G}$. Set
\[ n : =
\Bigg\lceil\frac{\theta_{\Op'\ps{G}}\big(X(A/F_{\infty})\big)+1}{e}\Bigg\rceil,\]
where $e$ is the ramification index of $\Op'/\Op$. Then for every
prime element $g$ of $\Op\ps{T}$ with $\pi^{n}|f-g$ such that
$\Op\ps{T}/g$ is isomorphic to $\Op'$, we have that
$X(\A[g]/F_{\infty})$ is torsion over $\Op'\ps{G}$, and that
$X(A/F_{\infty})(\pi')$ and $X(\A[g]/F_{\infty})(\pi')$ have the
same elementary representations. \ep

Note that by Lemma \ref{Greenberg equal complex}, this proposition
may be viewed as a refinement of \cite[Theorem 1.2(1), Corollary
4.37(1)]{B}). We now give the proof.

\bpf[Proof of Proposition \ref{big galois compare}]
 Let $g$ be a prime element of $\Op\ps{T}$ which satisfies the hypothesis
in the proposition. Let $\pi'$ be a prime element of $\Op'$ and
write $B= \A[g]$. It follows from Lemma \ref{big galois lemma} that
there is an isomorphism of $G_S(F)$-modules $A[\pi'^{en}]\cong
B[\pi'^{en}]$ which induces an isomorphism of
$\Gal(\bar{F}_v/F_v)$-modules $A_v[\pi'^{en}]\cong B_v[\pi'^{en}]$
for each prime $v$ of $F$ above $p$. By our hypothesis of $n$, we
have $en \geq \theta_{\Op'\ps{G}}(A)+1$. In particular, the
congruence hypothesis $\mathbf{(C_{\theta_{\Op'\ps{G}}(A)+1})}$
holds for $A$ and $\A[g]$. Hence the conclusion of the proposition
is now immediate from Theorem \ref{congruent}. \epf

We end the paper with a proposition which is immediate from an
application of Proposition \ref{MHG congruent}. This proposition is
a refinement of \cite[Proposition 8.6]{SS} and \cite[Corollary
4.37]{B} when the admissible $p$-adic Lie extension is of dimension
2.

\bp \label{big galois MHG}
 Let $F_{\infty}$ be an
admissible $p$-adic Lie extension of $F$, whose Galois group is a
pro-$p$ group of dimension 2 and has no elements of order $p$. Let
$f$ be a prime element of $\Op\ps{T}$ such that $\Op':=\Op\ps{T}/f$
is a maximal order. Set $A = \A[f]$ and suppose that
$X(A/F_{\infty})$ belongs to $\M_H(G)$. Set \[n : =
\Bigg\lceil\frac{\theta+1}{e}\Bigg\rceil,\] where $\theta =
\max\{\theta_{\Op'\ps{\Ga}}(A), \theta_{\Op'\ps{G}}(A)\}$ and $e$ is
the ramification index of $\Op'/\Op$. Then for every prime element
$g$ of $\Op\ps{T}$ with $\pi^{n}|f-g$ such that $\Op\ps{T}/g$ is
isomorphic to $\Op'$, we have that $X(\A[g]/F_{\infty})$ belongs to
$\M_H(G)$.
 \ep


\footnotesize
\begin{thebibliography}{99999}


\bibitem[AB]{AB} K. Ardakov and K. A. Brown, Primeness, semiprimeness and localisation in Iwasawa algebras,
\textit{Trans. Amer. Math. Soc.} 359(4) (2007) 1499--1515.

\bibitem[BZ]{BZ} T. Backhausz and G. Z\'{a}br\'{a}di,
Algebraic functional equations and completely faithful Selmer
groups, \textit{Int. J. Number Theory} \textbf{11} (2015)
1233--1257.

\bibitem[BS]{BS} R. Barman and A. Saikia, A note on Iwasawa
$\mu$-invariants of elliptic curves, \textit{Bull Braz Math Soc, New
Series} \textbf{41}(3) (2010) 399--407.

\bibitem[B]{B} P. Barth, Iwasawa theory for one-parameter families
of motives, \textit{Int. J. Number Theory} \textbf{9}(2) (2013)
257--319.

\bibitem[Bh]{Bh} A. Bhave, Comparison of the $\mu$-invariants of an abelian
variety and its dual abelian variety, arXiv:1305.3444v1 [math.NT].

\bibitem[Ch]{Ch} A. Chandrakant Sharma, Iwasawa invariants for the False-Tate
extension and congruences between modular forms, \textit{J. Number
Theory} \textbf{129} (2009) 1893--1911.

%\bibitem[CFKS]{CFKS} J.\ Coates, T.\ Fukaya, K.\ Kato and R.\ Sujatha Root numbers,
%Selmer groups and noncommutative Iwasawa theory, \textit{J.
%Algebraic Geom.} \textbf{19}(1) (2010) 19--97.

\bibitem[CF$^+$]{CFKSV} J.\ Coates, T.\ Fukaya, K.\ Kato, R.\ Sujatha
and O.\ Venjakob, The GL$_2$ main conjecture for elliptic curves
without complex multiplication, \textit{Publ. Math. IHES}
\textbf{101} (2005) 163--208.

\bibitem[CG]{CG} J. Coates and R. Greenberg, Kummer theory for abelian varieties over
local fields, \textit{Invent. Math.} \textbf{124} (1996) 129--174.

\bibitem[CS]{CS12} J. Coates and R. Sujatha, On the
$\mathfrak{M}_H(G)$-conjecture, in \textit{Non-abelian fundamental
groups and Iwasawa theory}, ed. J. Coates, M. Kim, F. Pop, M. Saidi
and P. Schneider, London Math. Soc. Lecture Note Ser. \textbf{393},
Cambridge Univ. Press, 2012, pp. 132--161.

\bibitem[CSW]{CSW} J. Coates, R. Sujatha and J -P Winterberger, On
the Euler-Poincar\'e characteristics of finite dimensional $p$-adic
Galois representations, \textit{Publ. Math. IHES} \textbf{93} (2001)
107--143.


\bibitem[CM]{CM} A. Cucuo and P. Monsky, Class numbers in
$\Zp^d$-extensions, \textit{Math. Ann.} \textbf{255}(2) (1981)
235--258.

\bibitem[DS+]{DSMS} J. Dixon, M. P. F. Du Sautoy, A. Mann and D. Segal, \textit{Analytic Pro-p Groups},
2nd edn, Cambridge Stud. Adv. Math. \textbf{38}, Cambridge Univ.
Press, Cambridge, UK, 1999.

\bibitem[EPW]{EPW} M. Emerton, R. Pollack and T. Weston,
Variation of Iwasawa invariants in Hida families, \textit{Invent.
Math.} \textbf{163} (2006) 523--580.

\bibitem[FK]{FK} T. Fukaya and K. Kato, A formulation of conjectures on $p$-adic zeta
functions in noncommutative Iwasawa theory, \textit{Amer. Math. Soc.
Transl. Ser. 2} \textbf{219}, 2006, 1--85.

\bibitem[GW]{GW} K. R. Goodearl and R. B. Warfield, \textit{An
introduction to non-commutative Noetherian rings}, London Math. Soc.
Stud. Texts \textbf{61}, Cambridge University Press, 2004.

\bibitem[Gr1]{G89} R. Greenberg, Iwasawa theory for $p$-adic representations, in
\textit{Algebraic Number Theory--in honor of K. Iwasawa}, ed. J.
Coates, R. Greenberg, B. Mazur and I. Satake, Adv. Std. in Pure
Math. \textbf{17}, 1989, pp. 97--137.

\bibitem[Gr2]{Gr94} R. Greenberg, Iwasawa theory for $p$-adic deformations
of motives, \textit{Proc. Sympos. Pure Math.} \textbf{55} (Part 2)
(1994) 193--223.

\bibitem[GV]{GV} R. Greenberg and V. Vatsal, On the Iwasawa invariants of elliptic
curves, \textit{Invent. Math.} \textbf{142} (2000) 17--63.

\bibitem[Ha]{Ha} Y. Hachimori, Iwasawa $\lambda$-invariants and
congruence of Galois representations, \textit{J. Ramanujan Math.
Soc.} \textbf{26}(2) (2011) 203--217.

\bibitem[Har1]{Har} M. Harris, $p$-adic representations arising from descent on abelian
varieties, \textit{Comp. Math.} \textbf{39} (1979) 177--245.

\bibitem[Har2]{Har2} M. Harris, Correction to $p$-adic
representations arising from descent on abelian varieties,
\textit{Comp. Math.} \textbf{121} (2000) 105--108

\bibitem[Ho1]{Ho} S. Howson, Euler characteristic as invariants of
Iwasawa modules, \textit{Proc. London Math. Soc.} \textbf{85}(3)
(2002) 634--658.

\bibitem[Ho2]{Ho2} S. Howson, Structure of central torsion Iwasawa
modules, \textit{Bull. Soc. Math. France} \textbf{130}(4) (2002)
507--535.

\bibitem[Hs]{Hs} M.-L. Hsieh, The algebraic functional equation of Selmer groups for CM
fields, \textit{J. Number Theory} \textbf{130} (2010), 1914--1924.

\bibitem[Im]{Im} H. Imai, A remark on the rational points of abelian
varieties with values in cyclotomic $\Zp$-extensions, \textit{Proc.
Japan Acad.} \textbf{51} (1975) 12--16.

\bibitem[Iw]{Iw} K. Iwasawa, On $\Gamma$-extensions of algebraic number fields,
\textit{Bull.
Amer. Math. Soc.} \textbf{65} (1959) 183--226.

\bibitem[JP]{JP} S. Jha and A. Pal, Algebraic functional equation
for Hida family, \textit{Int. J. Number Theory} \textbf{10}(7)
(2014) 1649--1674.

\bibitem[KT]{KT} Y. Kubo and Y. Taguchi, A generalization of a theorem of Imai
and its applications to Iwasawa theory, \textit{Math. Z.}
\textbf{275}(3-4) (2013) 1181--1195.

\bibitem[LLTT]{LLTT} K. F. Lai, I. Longhi, K. -S. Tan and F. Trihan,
Pontryagin duality for Iwasawa modules and abelian varieties,
arXiv:1406.5815 [math.NT].

\bibitem[Lam]{Lam} T. Y. Lam, \textit{Lectures on Modules and Rings},
Grad. Texts in Math. \textbf{189}, Springer, 1999.

\bibitem[Laz]{Laz} M.\ Lazard, Groups analytiques $p$-adiques, \textit{Pub. Math. IHES}
$\mathbf{26}$ (1965) 389-603.

\bibitem[Lim]{LimMu} M. F. Lim, Comparing the Selmer group of a $p$-adic
representation and the Selmer group of the Tate dual of the
representation, unpublished note, available at arXiv:1405.5289 [math.NT].

\bibitem[Mat]{Mat} K. Matsuno, Finite $\La$-submodules of Selmer groups of abelian
varieties over cyclotomic $\Zp$-extensions, \textit{J. Number
Theory} \textbf{99}(2) (2003) 415--443.

\bibitem[Mon]{Mon} P. Monsky, Fine estimate for the growth of $e_n$ in $\Zp^d$-extensions,
in \textit{Algebraic Number Theory--in honor of K. Iwasawa}, ed. J.
Coates, R. Greenberg, B. Mazur and I. Satake, Adv. Std. in Pure
Math. \textbf{17}, 1989, pp. 309--330.

\bibitem[Nek]{Nek} J. Nekov\'a\v{r}, \textit{Selmer Complexes},
Ast\'erisque \textbf{310}, 2006.

\bibitem[NSW]{NSW} J.\ Neukirch, A.\ Schmidt and K.\ Wingberg,
\textit{Cohomology of Number Fields}, 2nd edn, Grundlehren Math.
Wiss. \textbf{323}, Springer 2008.

\bibitem[Neu]{Neu} A. Neumann, Completed group algebras without zero divisors,
\textit{Arch. Math.} \textbf{51}(6) (1988) 496--499.

\bibitem[Oc]{Oc} T. Ochiai, On the two-variable Iwasawa main
conjecture, \textit{Comp. Math.} \textbf{142} (2006) 1157--1200.

\bibitem[Ser]{Ser} J -P. Serre, Sur la dimension cohomologique des groupes profinis,
\textit{Topology} \textbf{3} (1965) 413--420.

\bibitem[SS]{SS} S. Shekhar and R. Sujatha, On the structure of
Selmer groups of $\La$-adic deformations over $p$-adic Lie
extensions, \textit{Doc. Math.} \textbf{17} (2012) 573--606.

\bibitem[Su]{Su} R. Sujatha, Iwasawa theory and modular forms,
\textit{Pure Appl. Math. Q.} \textbf{2}(2) (2006) 519--538.

\bibitem[V1]{V02} O. Venjakob, On the structure theory of the Iwasawa algebra of a
$p$-adic Lie group, \textit{J. Eur. Math. Soc.} \textbf{4}(3) (2002)
271--311.

\bibitem[V2]{V03} O. Venjakob, A non-commutative Weierstrass preparation theorem
and applications to Iwasawa theory, \textit{J. reine angew. Math.}
\textbf{559} (2003) 153--191.

\bibitem[We]{We} T. Weston, Iwasawa invariants of Galois
deformations, \textit{Manuscripta Math.} \textbf{118}(2) (2005)
161--180.

\bibitem[Z1]{Z08} G. Z\'{a}br\'{a}di, Characteristic elements, pairings and functional
equations over the false Tate curve extension, \textit{Math. Proc.
Camb. Phil. Soc.} \textbf{144} (2008) 535--574.

\bibitem[Z2]{Z10} G. Z\'{a}br\'{a}di, Pairings and functional equations over the
GL$_2$-extension, \textit{Proc. London Math. Soc.} \textbf{101}(3)
(2010) 893--930.

\end{thebibliography}
\end{document}